\documentclass[preprint,12pt]{elsarticle}

\usepackage{amssymb}

\usepackage{xspace}
\usepackage{subfigure}
\usepackage{amsmath}
\usepackage{enumerate}
\usepackage{framed}
\usepackage{graphicx}
\usepackage{lscape}
\usepackage{bm}
\usepackage{color}
\usepackage{multicol}
\usepackage{amsthm}
\usepackage{amsthm}

\newcommand{\figref}[1]{{Figure~\ref{#1}}}

\newcommand{\secref}[1]{{Section~\ref{#1}}}
\renewcommand{\eqref}[1]{{(\ref{#1})}}

\newtheorem{theorem}{Theorem}[section]

\newtheorem{remark}[theorem]{Remark}

\begin{document}

\begin{frontmatter}



\title{A Stochastic Delay Model for Pricing Debt and Equity: Numerical Techniques and Applications.}

\author{Elisabeth Kemajou}
\ead{isakema@umn.edu}
\address{Department of Mathematics, University of Minnesota, USA}
\author{Antoine Tambue\corref{cor1}}
\ead{Antoine.Tambue@math.uib.no}
\cortext[cor1]{Address all correspondence to this author }
\address{ Department of Mathematics, University of Bergen, P.O. Box 7800, N-5020 Bergen, Norway}
\author {Salah Mohammed}
\ead{salah@sfde.math.siu.edu}
\address{Department of Mathematics, Southern Illinois University Carbondale, USA}


\begin{abstract}
In the  accompanied paper \cite{ElisaANTO}, a  delayed nonlinear model  for  pricing  corporate liabilities was developed.
Using self-financed strategy and duplication we were able to derive two Random Partial Differential Equations (RPDEs) 
describing the evolution of  debt and equity values of the corporate in the last delay period interval. 
In this paper, we provide numerical techniques to solve 
our delayed nonlinear model along with the corresponding RPDEs modeling the debt and equity values of the corporate. 

Using financial  data from some firms, we compare  numerical solutions from both our nonlinear model and classical  Merton model \cite{GD2} to the real corporate data.
From this comparison, it comes up that in corporate  finance  the past dependence of
the firm value process may be an important feature and therefore  should not be ignored.

\end{abstract}


\end{frontmatter}
\section{Introduction}
Due to the remarkable growth of the credit derivatives market, the interest in corporate
claim value models and risk structure has recently increased.
Financial distress tends to be an important factor in many corporate decisions. The two
main sources of financial distress are corporate illiquidity and insolvency. In his paper \cite{Gryglewicz},
Gryglewicz explains how changes in solvency affect liquidity and also how liquidity concerns
affect solvency via capital structure choice. Corporate solvency is the ability to cover debt
obligations in the long run. Uncertainty about average future profitability, with financial
leverage, generates solvency concerns. Corporate insolvency may lead to corporate reorganization 
or to bankruptcy of the firm in the worst case.
Corporate bankruptcy is central to the theory of the firm. A firm is generally considered
bankrupt when it cannot meet a current payment on a debt obligation. In this event, 
the equity holders lose all claims on the firm, and the remaining loss which is the difference between
the face value of the fixed claims and the market value of the firm, is supported by the debt
holders. 
In the literature of corporate finance,  Merton \cite{GD2}
appears to be the main pioneers in the derivation of formulas for corporate claims. This model is a dual  of Black and Scholes model \cite{GD1} for stock price.
 Merton \cite{GD2} further analyzed the risk structure of interest rates. More specifically, he found the relation
between corporate bond spreads and government bond, and attempted to determine a valid
measure of risk. He also developed the deterministc partial differential equation  modelling the debt and equity of the firm.  
The assumption of constant volatility in  the original Black-Scholes  and Merton models from which most claims derivations are inspired, is
incompatible with derivatives prices observed in the market (see \cite{Bernard,BD, Blattberg,scott} and the references therein). For stock price, two alternative theories are mostly used to overcome the constant volatility drawback.
The first approach  sometime called  {\tt level-dependent volatility} describes the stock price as a diffusion with level dependent volatility \cite{Bensoussan}.
The second approach  sometime called  {\tt stochastic volatility} defines the volatility as an autonomous diffusion driven by a second Brownian motion.
 \footnote{In  the sense that the first Brownian motion drives the asset price}

In \cite{firstdelay}, a new class of   nonconstant volatility model  which can be extended to include the first of the above approaches,
that  we called  { \tt delayed model } is  introduced and further study in \cite {GD4,stra} for options prices.  This model shows that  the past dependence of
the stock price process is an important feature and therefore  should not be ignored.  The main  goal of this model is to make volatility self--reinforcing.
Since the volatility is defined in terms of past behavior of the asset price,  the self--reinforcing  is high, precisely when there have
been large movements in the recent past (see \cite{firstdelay}). This is designed to reflect real--world perceptions of
market volatility, particularly if practitioners are to compare historic volatility with implied.

 Following the duality  between   the  stock price \cite{GD1} and corporate finance \cite{GD2}, we have recently  introduced 
 in  \cite{ElisaThesis,ElisaANTO} the nonlinear delayed model in debt and guarantee.  Using self-financed strategy and replication 
we  established  that debt  value and equity value
 follow  two similar Random Partial differential Equations  (RPDEs) within  the last delay period interval \footnote{The final  
time interval with length equal to the time delay}.
 The analytical solution  of our nonlinear model and  RPDEs are unknown in general case and therefore   numerical techniques  are needed.

In recent years, the computational complexity of mathematical models employed in financial mathematics has witnessed a tremendous growth 
(see \cite{P3,option,optionbound} and references therein).
The aims of this paper is to solve numerically 
   our delayed nonlinear model for firm market value along with the corresponding  RPDEs, using real data from firms.  
   Comparison will be  done with classical Merton model.  To the  best of our knowledge such comparison has  not yet been done in the financial literature.
   Two major comparisons will be preformed:  the market value of each corporate and its equity value (or its debt value).
   We will first   approximate the volatility of each  corporate, afterward solve  numerically our nonlinear model  for the market value of 
the corporate along with the corresponding  Merton model  using the $\theta-$ semi  implicit Euler Maruyama
scheme  to obtain sample numerical solutions.  Monte carlo method will be thereafter used to approximate the  mean numerical solution of each model.
The meam numerical value  from each model (our nonlinear model and Merton model) will be therefore compared with the real market value ($V$) of the corporate. 
For  debt value (or equity value) solutions of  RPDEs established in the accompanied paper \cite{ElisaANTO}, efficient numerical scheme based on finite volume-finite 
difference methods (discretization  respect to the firm value $V$) and exponential integrator (discretization respect to the time $t$) will be used.  
Recently, exponential integrators have been used  efficiency in many applications in  porous media flow \cite{Antoine,TLGspe1,ATthesis,Antoine2012,Antoine2012b},  but are not yet well spread in finance.  
The same numerical technique is also used to solve deterministic Partial Differential Equations (PDEs) modeling debt  value or equity value in Merton model. 
Comparisons are done with the real data from  firms  for each model (our delay model and Merton model).

From our comparison, it comes up that in corporate  finance  the past dependence of
the firm value process is an important feature and therefore  should not be ignored. The main goal of this paper is to call for further attention into the possibility 
of modeling  market  value of the firm  with nonlinear delayed stochastic differential equations. 

The paper is organized as follows. In \secref{model}, we recall our delayed nonlinear model for corporate claims as presented in \cite{ElisaANTO} 
along with the Merton model \cite{GD2}.  In \secref{numeric}, numerical techniques for our  delayed nonlinear model are provided. We first present 
the  $\theta-$ semi implicit Euler Maruyama for  the firm market value  $V$ and provide numerical experimentations for both our nonlinear 
model and Merton model using real data for some firms. We end this section by providing  numerical technique to solve efficiently our 
(RPDEs) modeling the debt and equity of the firm along with numerical experimentations for the two models (our delayed nonlinear model and Merton model) 
with real data  for some firms. The conclusion is provided in \secref{conclusion}.



\section{Stochastic delay model for corporate claims}
 \label{model}
Here we present the  stochastic delay model  formulated  in  the accompanied paper \cite{ElisaANTO} along with  Random Partial Differential Equation (RPDE) that should satisfy  any claim.
We assume that:
\begin{itemize}
 \item[$A_1$] The value of the company is unaffected by how it is financed (the capital structure irrelevance principle).
 \item[$A_2$]  The market value of firm at time $t \in [0,T]$, $V(t)$, follows a nonlinear Stochastic Delay
 Differential Equation (SDDE)
\begin{eqnarray}
\label{model}
\left\lbrace \begin{array}{l}
 dV(t)=(\alpha V(t)V(t-L)-C)dt+g(V(t-L))V(t)dW(t) \\
 \newline\\
 V(t)=\varphi(t),\,\, t\in [-L,0]\\
\end{array} \right.
\end{eqnarray}
on a probability space $(\Omega,\mathcal{F},P))$ with a filtration $(\mathcal{F}_{0\leq t\leq T})$ satisfying the usual conditions.
\end{itemize}
 where  $\alpha$ is the  riskless interest rate of return on the firm per unit time, $C$ 
is the total amount payout by the firm per unit time to either the shareholders or claims-holders (e.g. dividends or interest payments) 
if positive, and it is the net amount received by the firm from new financing if negative.  
The constant $L$ represents the past length while $T$ is the maturity date.
The function  $g:\mathbb{R}\rightarrow\mathbb{R}$ 
is a continuous representing the volatility function on the firm value per unit time.
The initial process $\varphi : \Omega \rightarrow C([-L,0],\mathbb{R})$ is $\mathcal{F}_0$-measurable with respect to the 
Borel $\sigma$-algebra of $C([-L,0],\mathbb{R})$, 
actually $\varphi$ is the past value of the firm.
The process $W$ is a one dimensional standard Brownian motion adapted to the filtration $(\mathcal{F}_t)_{0\leq t\leq T}$.

Notice that  $C$ and $\alpha $ can be time dependent functions, in which case they should be measurable and integrable in the interval $[0,T]$.

The results ensuring the feasiblility of   the price model  \eqref{model} is given in \cite{ GD2, ElisaThesis}.
 Following the work in \cite{GD2}, in order the RPDE which must be satisfied by any security whose value can be written as a function of the
value of the firm and time, we assume that any claim with market value $Y(t)$ (which can be replicated using self-financed strategy) at time  $t$ with $Y(t)=F(V(t),t)$ follows
a nonlinear stochastic delay differential equation
\begin{eqnarray}\label{eq1}
\left\lbrace \begin{array}{l}
  dY(t)=(\alpha_y Y(t) -C_y)dt+g_y(Y(t-L))Y(t)dW_y(t), \,\, t\in [0,T]\\
  \newline\\
  Y(t)=\varphi_y(t),\,\, t\in [-L,0],
  \end{array} \right.
\end{eqnarray}
on a probability space $(\Omega,\mathcal{F},P))$.
where 
$\alpha_y$ is the constant riskless interest rate of return per unit time on this claim; $C_y$ is the amount payout 
per unit time to this claim; $g_{y}: \mathbb{R}\rightarrow \mathbb{R}$ is a continuous function representing the volatility 
function of the return on this claim per unit time; the initial process $\varphi_y : \Omega \rightarrow C([-L,0],\mathbb{R})$ 
is $\mathcal{F}_0$-measurable with respect to the Borel $\sigma$-algebra of $C([-L,0],\mathbb{R})$. The functions $C_y(t)$ and $\alpha_y(t)$ 
are measurable and integrable in the interval $[0,T]$.
The process $W_y$ is a one dimensional standard Brownian motion adapted to the filtration $(\mathcal{F}_t)_{0\leq t\leq T}$.
For any claim
$Y(t)=F(V(t),t)$ where $F$ is twice continuously differentiable with respect to $V$ and once differentiable with respect to $t$, 
we have proved in \cite{ElisaANTO} that the following (RPDE) should be satisfied  
\begin{eqnarray}\label{eq:3008}
\dfrac{1}{2}g^2(V(t-L))v^2F_{vv}+(rv-C)F_v+F_t-rF+C_y=0,\\
 \qquad \qquad \,(t,v)\in [T-L,T] \times \mathbb{R}^{+}
\end{eqnarray}
with
$$F_t(v,t)=\frac{\partial F(v,t)}{\partial t},\, F_v(v,t)=\frac{\partial F(v,t)}{\partial v},\,F_{vv}(v,t)=\frac{\partial^2 F(v,t)}{\partial v^2}.$$

By setting 
\begin{eqnarray}
\label{en}
V(t) = F(V(t),t)+f(V(t),t),
 \end{eqnarray}
 where $f(V(t),t)$ is the value of the equity, $F(V(t),t)$ the value of debt a any time $t$  before the maturity and $r$ is the instantaneous
riskless rate of interest.
 We have obtained in \cite{ElisaANTO} the following  two final value problems for debt and equity, linked by 
 \eqref{en}
 \begin{eqnarray}
\label{eqoub}
\left \lbrace \begin{array}{l}
\dfrac{1}{2}g^2(V(t-L))v^2F_{vv}+(rv-C)F_v+F_t-rF+C_y =0,\,\,  t \in (T-L,T)\\
\newline\\
F(v,T)=\min[V,B],\,\,\,\, v>0\label{eqou1b}\\
\newline\\
F(0,t)=0, \,\,\,\, F(v,t)\sim Be^{-r(T-t)}, \text{ as }v\rightarrow\infty, \label{eqou2b}
\end{array} \right.
\end{eqnarray}
and
\begin{eqnarray}
\label{eqou}
\left \lbrace \begin{array}{l}
\dfrac{1}{2}g^2(V(t-L))v^2f_{vv}+(rv-C)f_v -rf +C-C_y+f_t=0,\,\, t \in ( T-L,T)\\
\newline\\
f(v,T)=\max(v-B,0),\,\,\,\,\, v>0 \label{eqou1}\\
\newline\\
f(0,t)=0, \,\,\,\, f(v,t)\sim v-Be^{-r(T-t)}, \text{ as } v\rightarrow\infty, \label{eqou2}
\end{array} \right.
\end{eqnarray}
where $B$ is   the promised value the firm must pay to the debtholders at the maturity date $T$.

\begin{remark}
 The classical Merton model \cite{GD2} assumes that  the value of the firm at time $t \in [0,T]$, $V(t)$, follows  a  Stochastic 
 Differential Equation (SDE)
\begin{eqnarray}
\label{modelM}
\left\lbrace \begin{array}{l}
 dV=(\alpha V-C)dt+ \sigma V dW \\
 \newline\\
 V(0) \,\, \text{given},\\
\end{array} \right.
\end{eqnarray}
where $\sigma^2$ is the  constant instantaneous variance of the
return on the firm per unit time. In  this case, the equity  value $f$ should satisfy  the following  deterministc PDE

\begin{eqnarray}
\label{eqouM}
\left \lbrace \begin{array}{l}
\dfrac{1}{2}\sigma^2v^2f_{vv}+(rv-C)f_v -rf +C-C_y+f_t=0,\,\, t \in ( 0,T)\\
\newline\\
f(v,T)=\max(v-B,0),\,\,\,\,\, v>0 \label{eqou1M}\\
\newline\\
f(0,t)=0, \,\,\,\, f(v,t)\sim v-Be^{-r(T-t)}, \text{ as } v\rightarrow\infty. \label{eqou2M}
\end{array} \right.
\end{eqnarray}

\end{remark}

\section{Numerical techniques and applications}
\label{numeric}
 \subsection{Presentation of the data set and volatility estimation}
%

The data on stock returns come from the  Center for  Research in Securty Prices(CRSP) database: http://www.crsp.com/ while
those on debt values are from the Research Insight/Compustat database (http://www.compustat.com/). 
More data include firms that had valid data for all 20 years from 1991-2010 and including:
\begin{enumerate}
 \item  The risk free rate $r$ per year,  which is the average monthly yield on US T-Bills for that year (the same for all firms each year).
\item  The standard deviation of daily returns $\sigma$ per  year for each firm. 
\item  The number of daily returns $N$ used to compute $\sigma$ for each firm each year (this is set to be at least 150).
\item The total book value of debt $B$ (in 1,000,000’s).   
\item The total value of the firm’s assets $V$  (in 1,000,000’s).
\item The total amount  $C$ (in 1,000,000’s) payout by the firm per unit time to either the shareholders or claims-holders for 10 years (2000-2010).
\item The  total amount  $C_y$  (in 1,000,000’s) payout  per unit time  for  the debt  within 10 years (2000-2010).
\end{enumerate}
In  fact  the data set we have used  include  all the parameters  that we need  
to solve either the stochastic differential equations \eqref{model} \& \eqref{modelM}, or the RPDE \eqref{eqoub} \& \eqref{eqou} and the PDE \eqref{eqouM}.

 All the simulation is performed in Matlab 7.7.
In most of our simulations, the data between 1991-2000.5  are used as memory data while those between  2000.5-2010 are used  as 
 the future data i.e.  the data that we want our model to predict. 

 To estimate the volatility function $g$,  we  use  the quadratic or  linear interpolation of the memory part of data $\sigma$.
 As in \cite{Bernard}, the quadratic form of the volatility is motivated by the fact that the implied volatility in Black-Scholes model has a parabolic shape.
 The volatility function  $g$  can also be estimateed by  using  the splines interpolation of the memory part of the data $ \sigma$.

As we only have yearly data set, we use also the interpolation to have more data  set if need as the numerical schemes usually need
small time step (then more data set) to ensure their stabilities.

\subsection{ Numerical approximation of the corporate market value}
\subsubsection{ The  $\theta-$ semi implicit Euler-Maruyama scheme}
Here we consider  the stochastic equations \eqref{model} and \eqref{modelM} within the time interval $[0,T]$, where the higher value of $T$ is 9.5 corresponding to the year $2010$. 
The time unit being the year. Indeed the values of $\alpha$ and $C$ are time depending, we therefore
consider those values as  two time depending functions, which are constant within  each year interval.
The goal here is to use the mean numerical solutions  of the firm $V$ as the forecasting values of firm in the interval $[2000.5,2000.5+T]$.
As  the real firm value of the companies are already known in that interval, the aim is to see how close are the forecasting firm values (from  numerical methods)
comparing to the real firm values  from financial industries. Recall that our nonlinear model  used the memory data within the interval $[1991, 2000.5]$.
We solve  numerically our nonlinear model \eqref{model}  for the value of 
the company and Merton model \eqref{modelM} \footnote{ For  constant $\alpha$ and $C$ the exact solution is well known as this is the same as Black Sholes model for stock price}  
in time  using the $\theta-$ implicit Euler Maruyama
scheme in order to obtain numerical sample solutions.  Monte carlo method is thereafter used to approximate the  mean numerical solution of each model.
The mean numerical value  from each model  will be therefore compared with the real company value $V$. 

  The  $\theta-$ semi implicit Euler-Maruyama scheme applied to  \eqref{model} is given by
\begin{eqnarray}
\label{scheme}
 V_{n+1}&=&V_{n}+\Delta T\left[\theta(\alpha_{n+1}V_{n+1}V_{n-m+1}-C_{n+1})+(1-\theta)(\alpha_{n}V_{n}V_{n-m}-C_{n})\right] \nonumber \\
\newline
    && \;\;+ g(V_{n-m}) V_{n} \Delta W_{n}\;\;\;\; n=1 .......,M,\;\; 0 \leq \theta \leq 1,\;\;\; L=m \Delta T. 
\end{eqnarray}
where  $\Delta T= T/M$ is the time step size,  $M$ the total number of time subdivision,  $V_{n}$ 
is the approximation of $V(t_{n}),\,t_{n}=n\Delta T$, $\alpha_{n}=\alpha(t_{n})$, $C_{n}=C(t_{n})$, and 
$$\Delta W_{n}= W(t_{n+1})- W(t_{n})$$  are standard Brownian increments,
independent 
identically distributed \\$ \sqrt{\Delta T}  \mathcal{N}(0,1)$ 
random variables.
For $\theta=0$, we have the classical Euler-Maruyama scheme which is less  numerical stable than the semi implicit Euler-Maruyama with 
$\theta=1$, that we will use in our simulation. To ensure the convergence of the numerical \eqref{scheme} toward the unique solution of  \eqref{model}, 
the volatility function
$g$ need to be globally Lipschitz, or localy Lipschitz and bounded \cite{stra}.
 These conditions are sufficient conditions for the convergence and not necessary conditions 
since the scheme can converge for some functions not verifying these conditions.

To approximate the expected value (mean) of  the process $V$, we use the  Monte Carlo method to compute the mean 
of the numerical  samples  from \eqref{scheme}. The  Monte Carlo method can also be used to approximate any moment of the process $V$.
 \subsubsection{Application with corporate data}
\begin{figure}[H]
 \subfigure[]{
   \label{FIG01a}
   \includegraphics[width=0.5\textwidth]{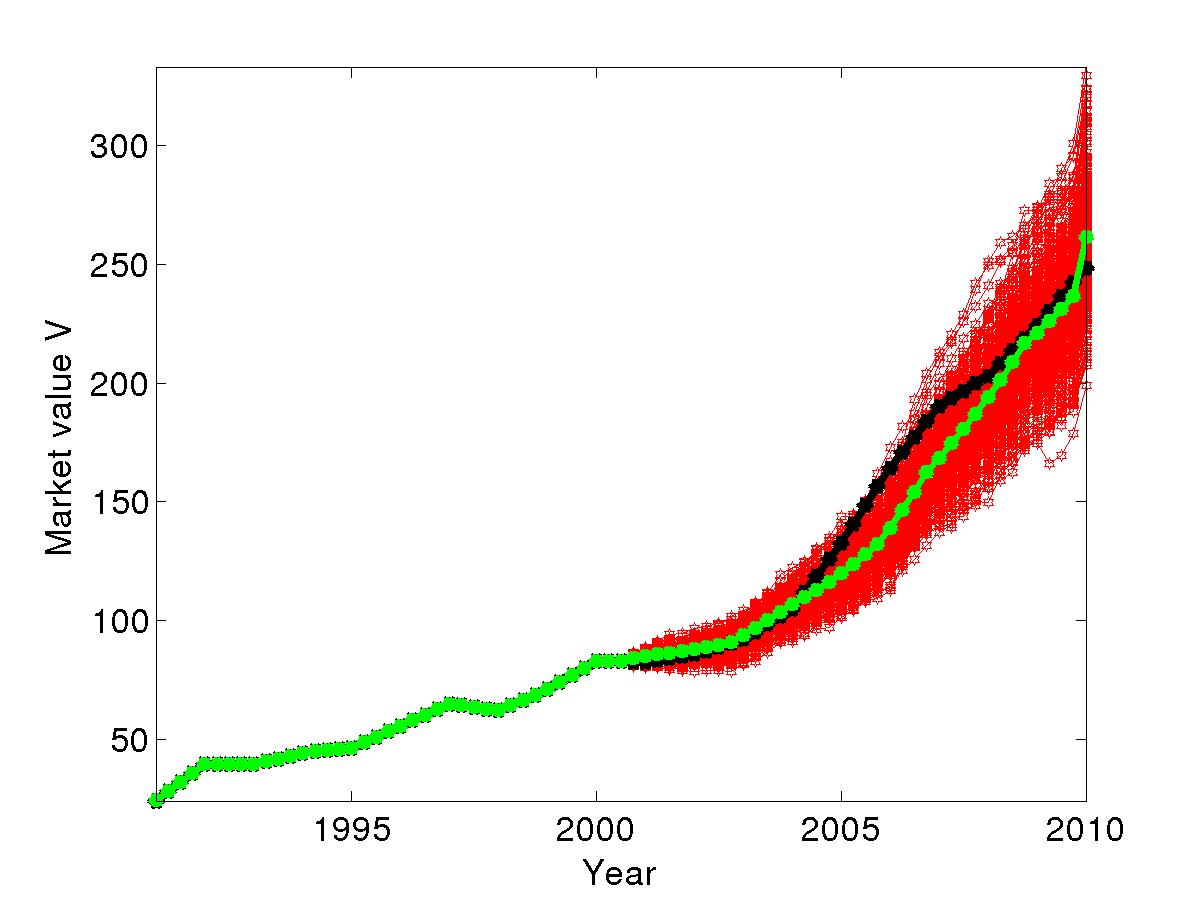}}
   \hskip 0.01\textwidth
   \subfigure[]{
   \label{FIG01b}
   \includegraphics[width=0.5\textwidth]{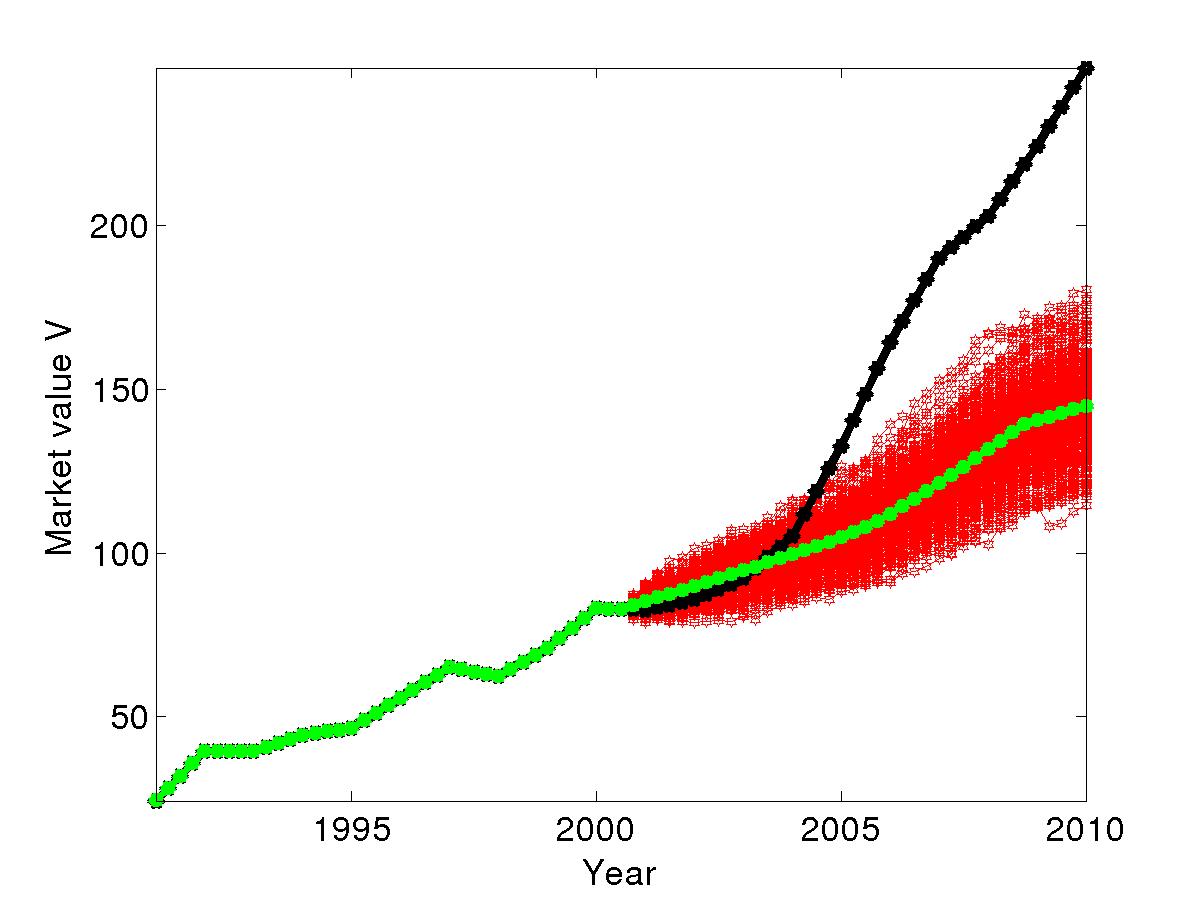}}
   \hskip 0.01\textwidth
   \subfigure[]{
   \label{FIG01c}
   \includegraphics[width=0.5\textwidth]{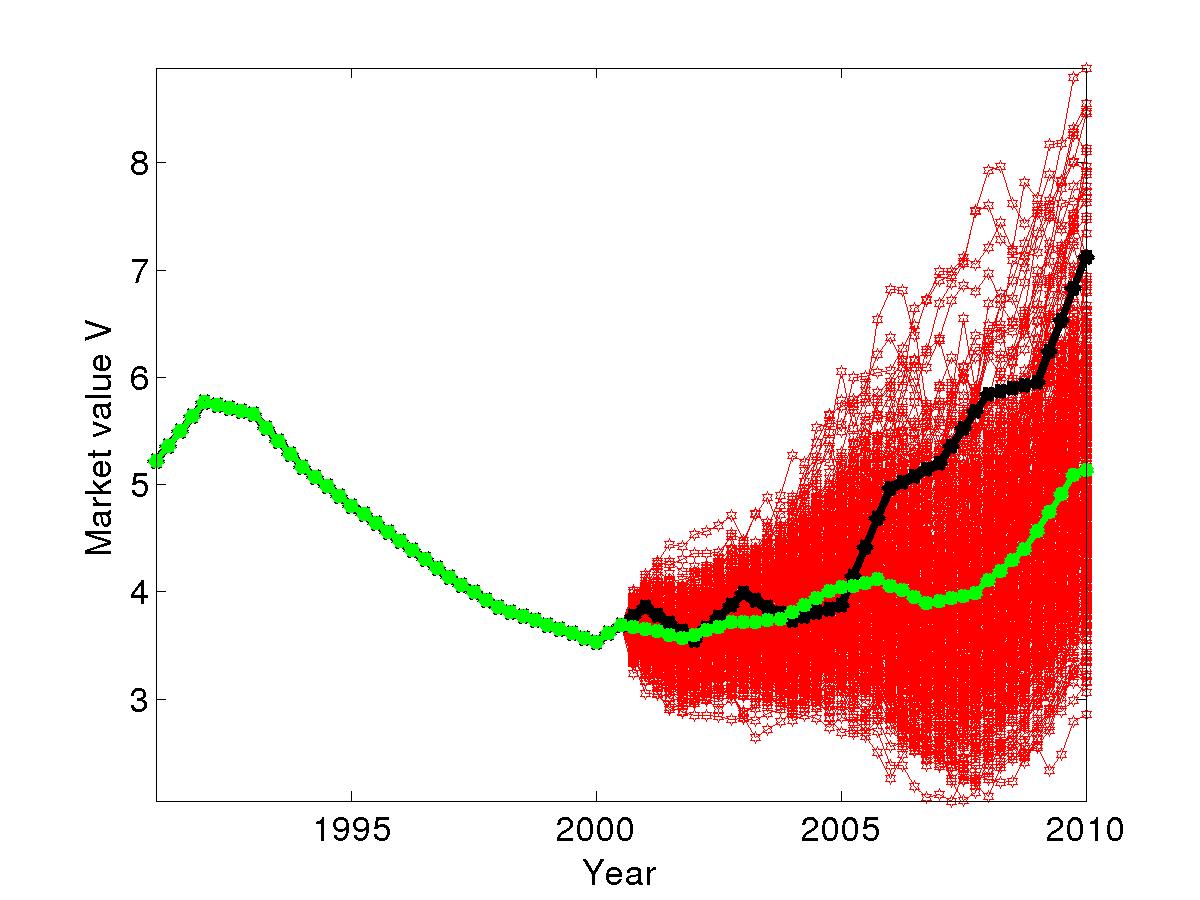}}
   \hskip 0.01\textwidth
   \subfigure[]{
   \label{FIG01d}
   \includegraphics[width=0.5\textwidth]{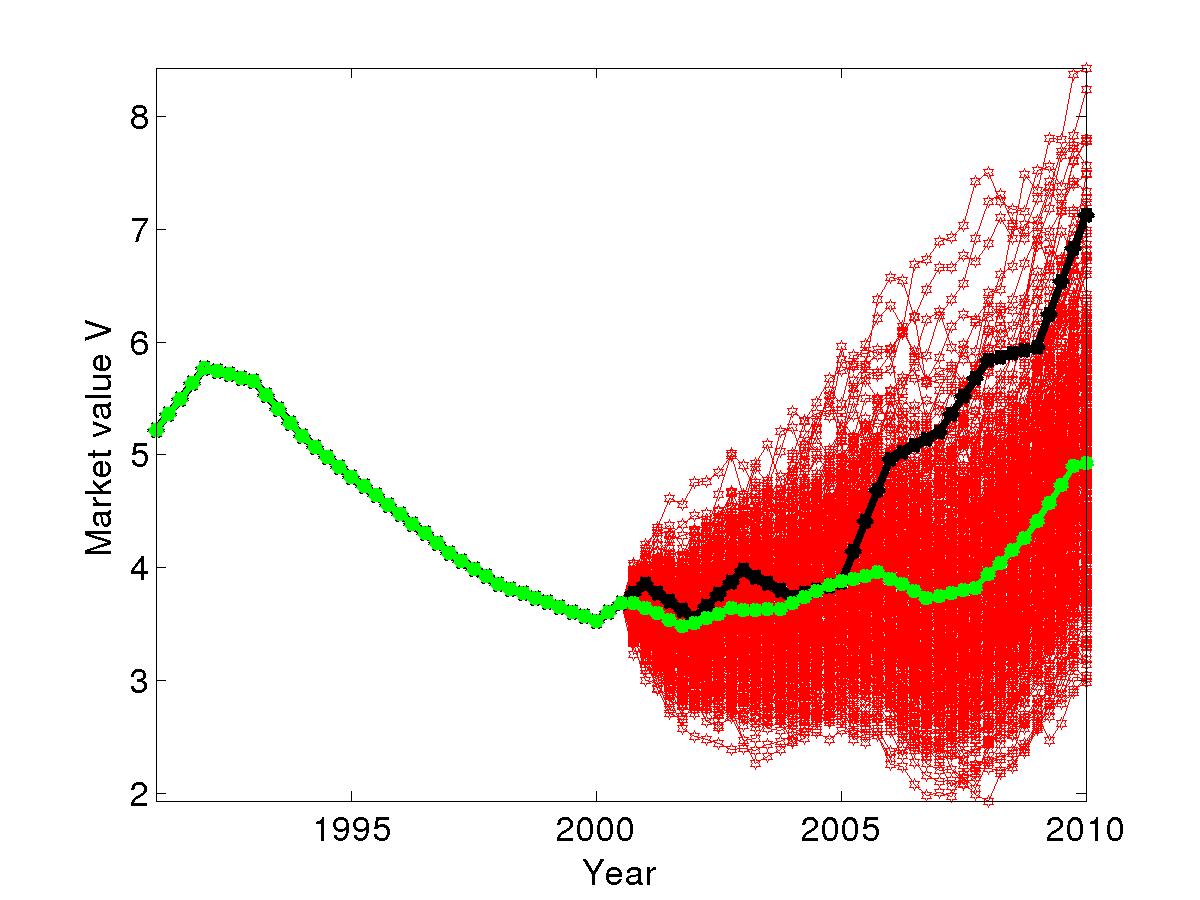}}
 \caption{The graphs at the left ((a),(c)  respectively for corporates $C_3$ and $C_6$ ) correspond to the delayed model  while the graphs at the right ((b),(d) 
 respectively for corporates $C_3$ and $C_6$) correspond to  Merton model.
    We aslso take $T=L=9.5$ and the function $g$ is the quadratic interpolation of the  standard deviation of daily returns $\sigma$ in the memory part. 
We  have plotted 400 samples of the numerical solution along with the expectation (the means)
 of the numerical solution (green curves). The curves of the real data of the firm  market value $V$ as a function of time  are in black (black thick curves).}
 \label{FIG01}
\end{figure}

\begin{figure}[H]
 \subfigure[]{
   \label{FIG011a}
   \includegraphics[width=0.5\textwidth]{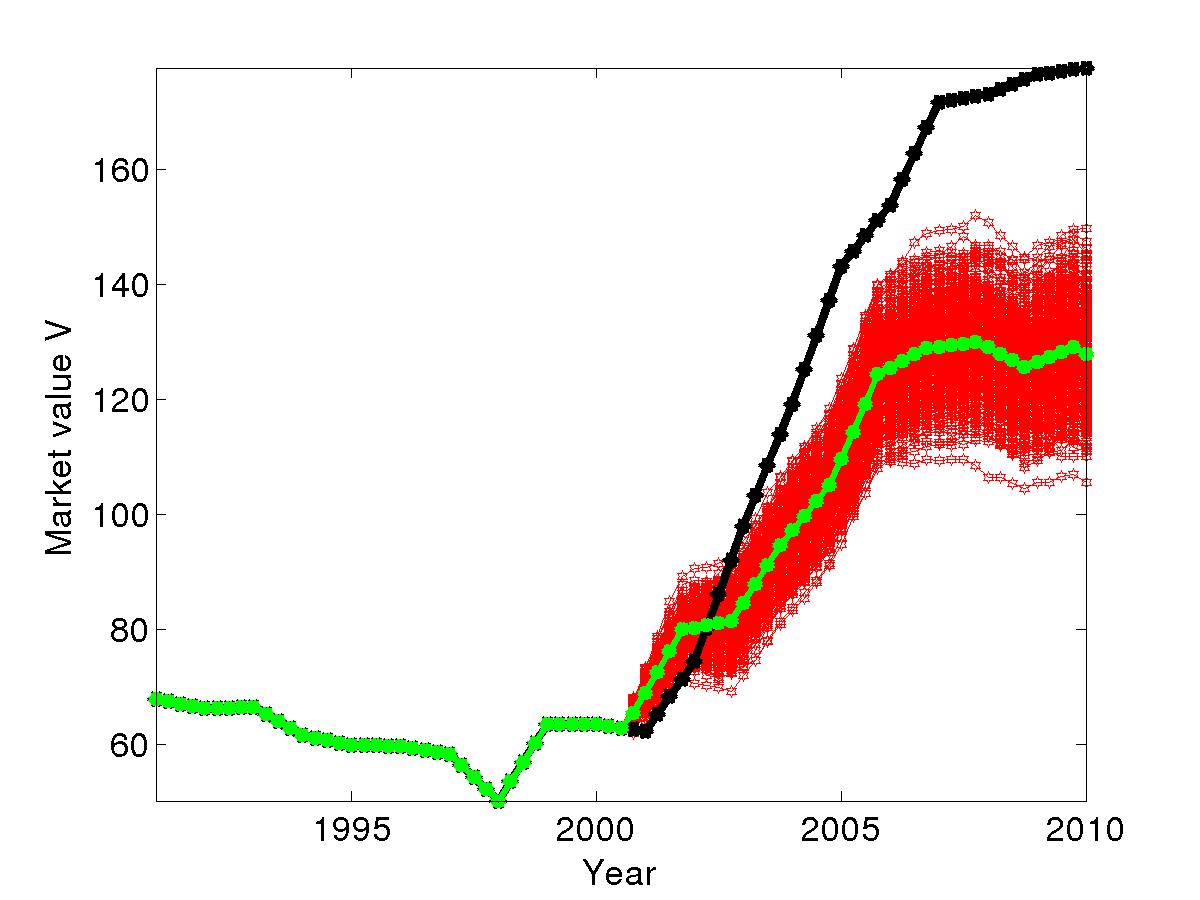}}
   \hskip 0.01\textwidth
   \subfigure[]{
   \label{FIG011b}
   \includegraphics[width=0.5\textwidth]{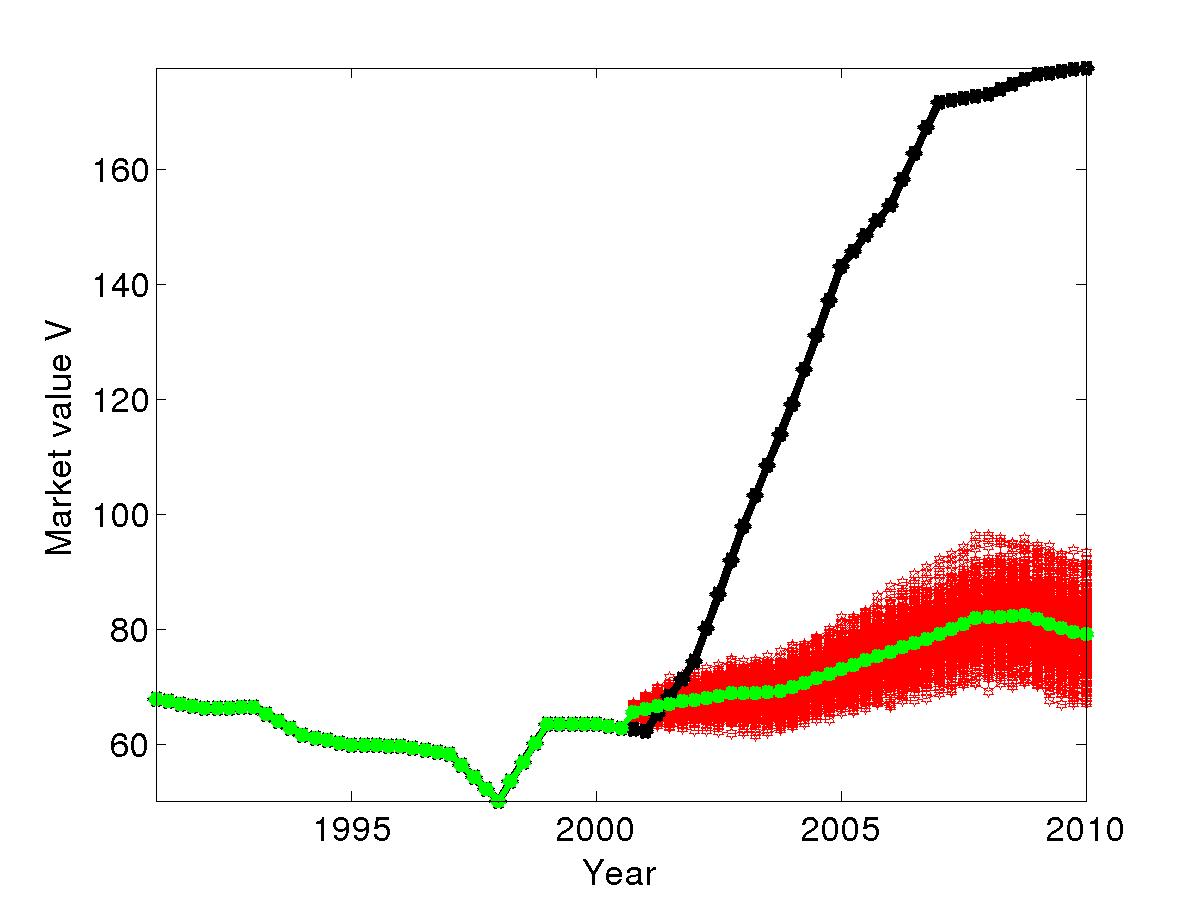}}
   \hskip 0.01\textwidth
   \subfigure[]{
   \label{FIG011c}
   \includegraphics[width=0.5\textwidth]{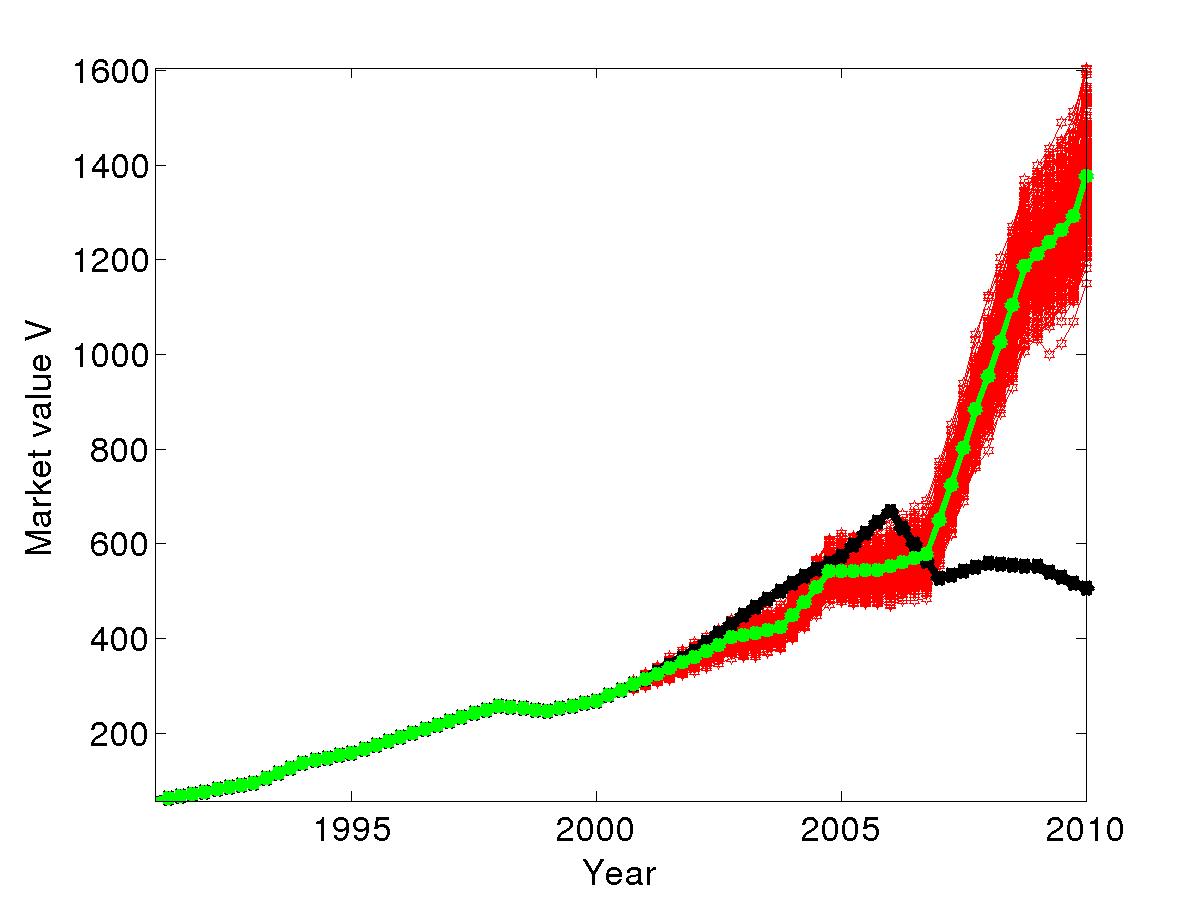}}
   \hskip 0.01\textwidth
   \subfigure[]{
   \label{FIG011d}
   \includegraphics[width=0.5\textwidth]{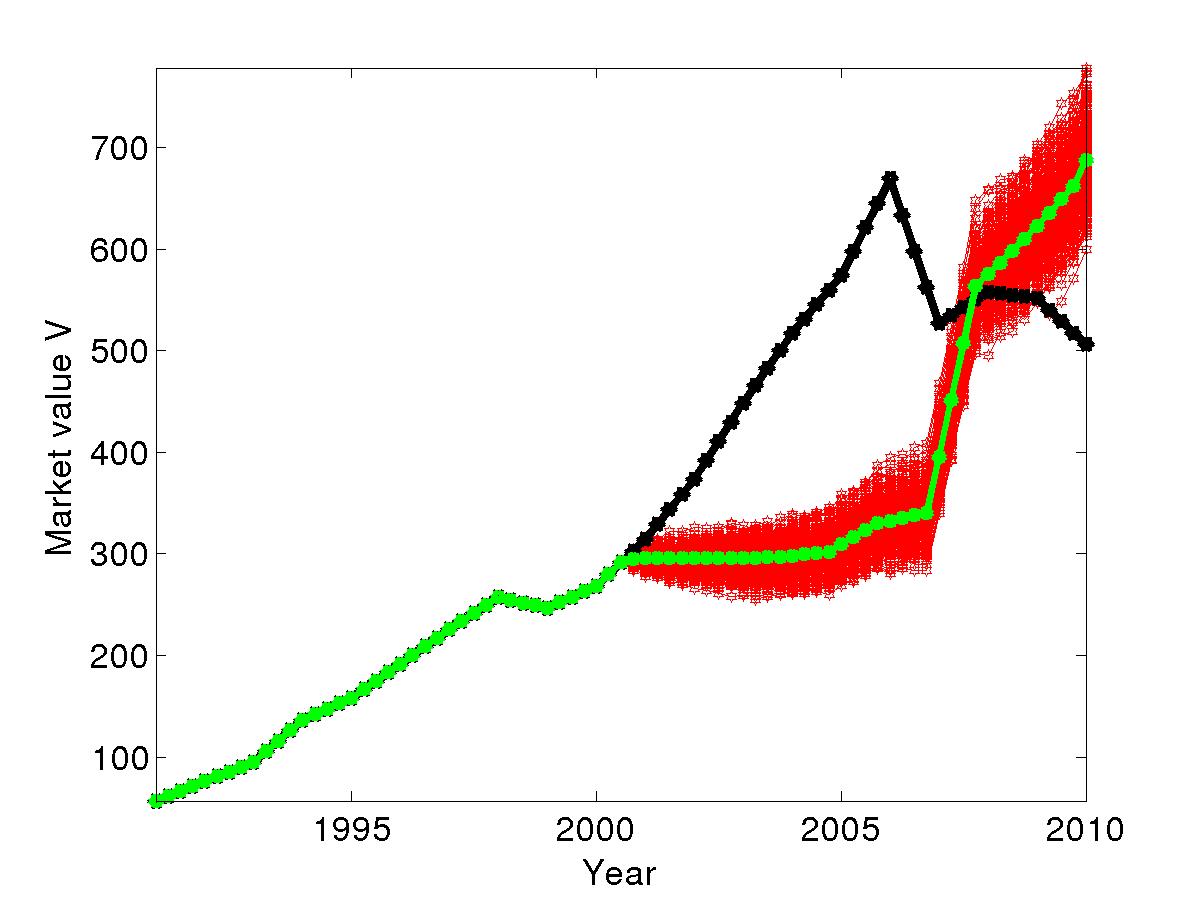}}
   \caption{ The graphs at the left ((a),(c)  respectively for corporates $C_5$ and $C_4$ ) correspond to the delayed model  while the graphs 
   at the right ((b),(d) respectively for corporates $C_5$ and $C_4$) correspond to  Merton model.
    We also take $T=L=9.5$ and the function $g$ is the quadratic interpolation of the  standard deviation of daily returns $\sigma$ in the memory part. 
We  have plotted 400 samples of the numerical solution along with the expectation (the means)
 of the numerical solution (green curves). The curves of the real data of the firm  market value $V$ as a function of time  are in black (black thick curves)}
   \label{FIG011}
\end{figure}

\begin{figure}
 \subfigure[]{
   \label{FIG02a}
   \includegraphics[width=0.5\textwidth]{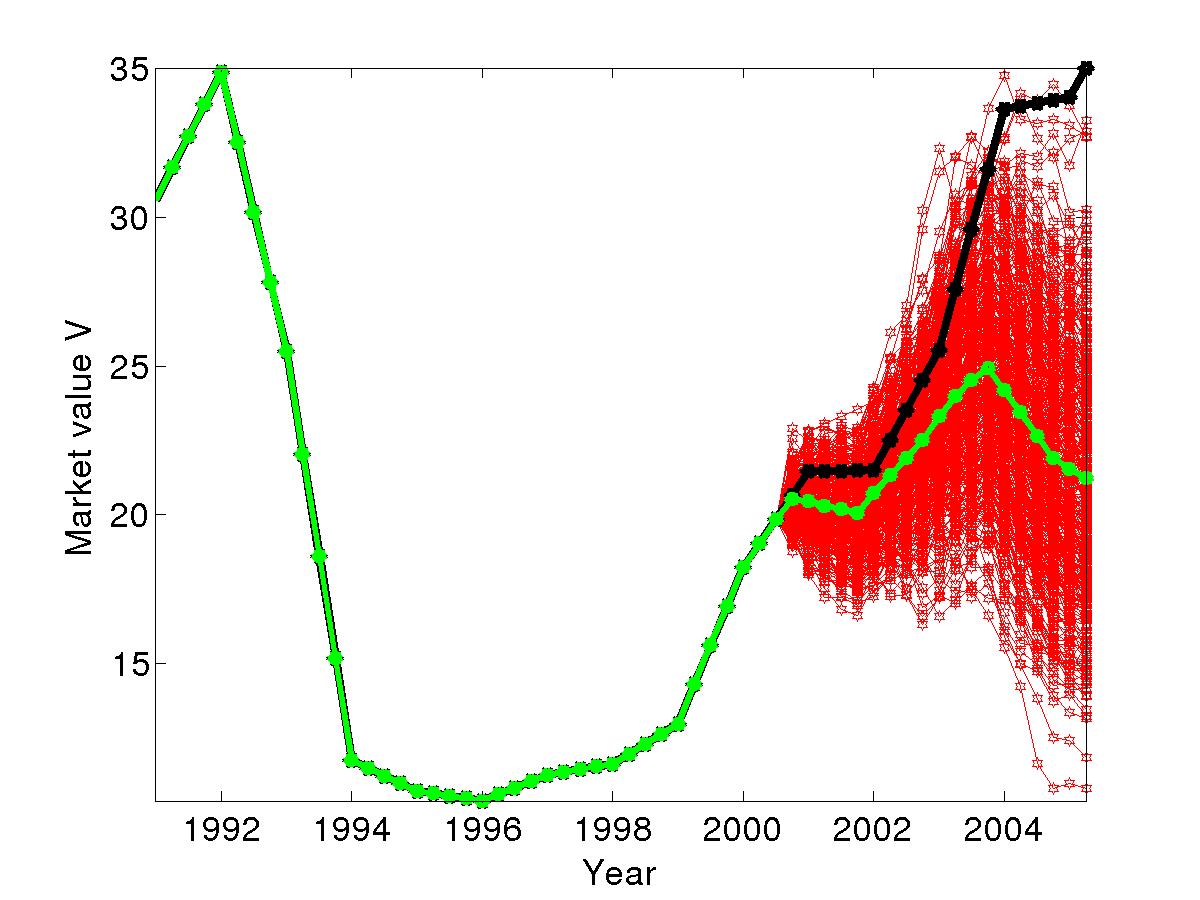}}
   \hskip 0.01\textwidth
   \subfigure[]{
   \label{FIG02b}
   \includegraphics[width=0.5\textwidth]{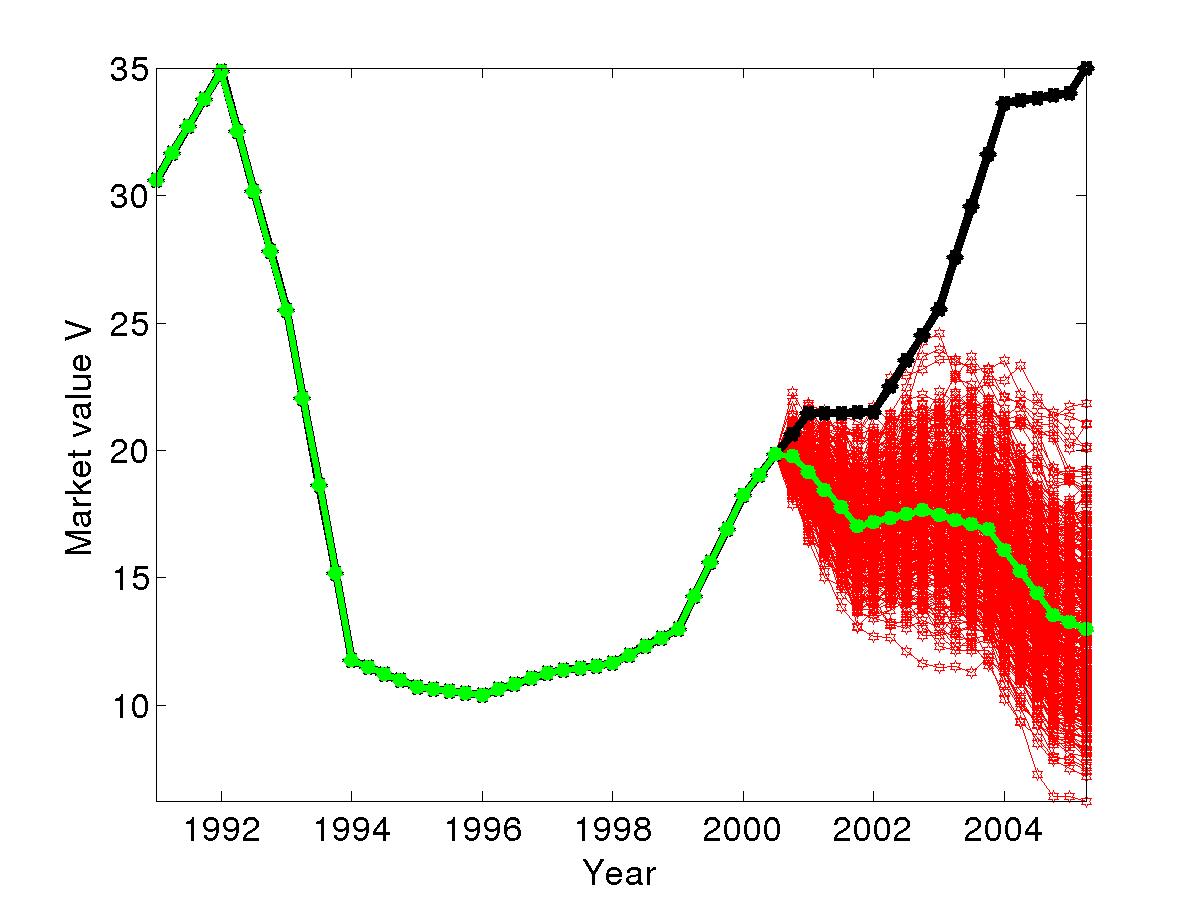}}
   \hskip 0.01\textwidth
   \subfigure[]{
   \label{FIG02c}
   \includegraphics[width=0.5\textwidth]{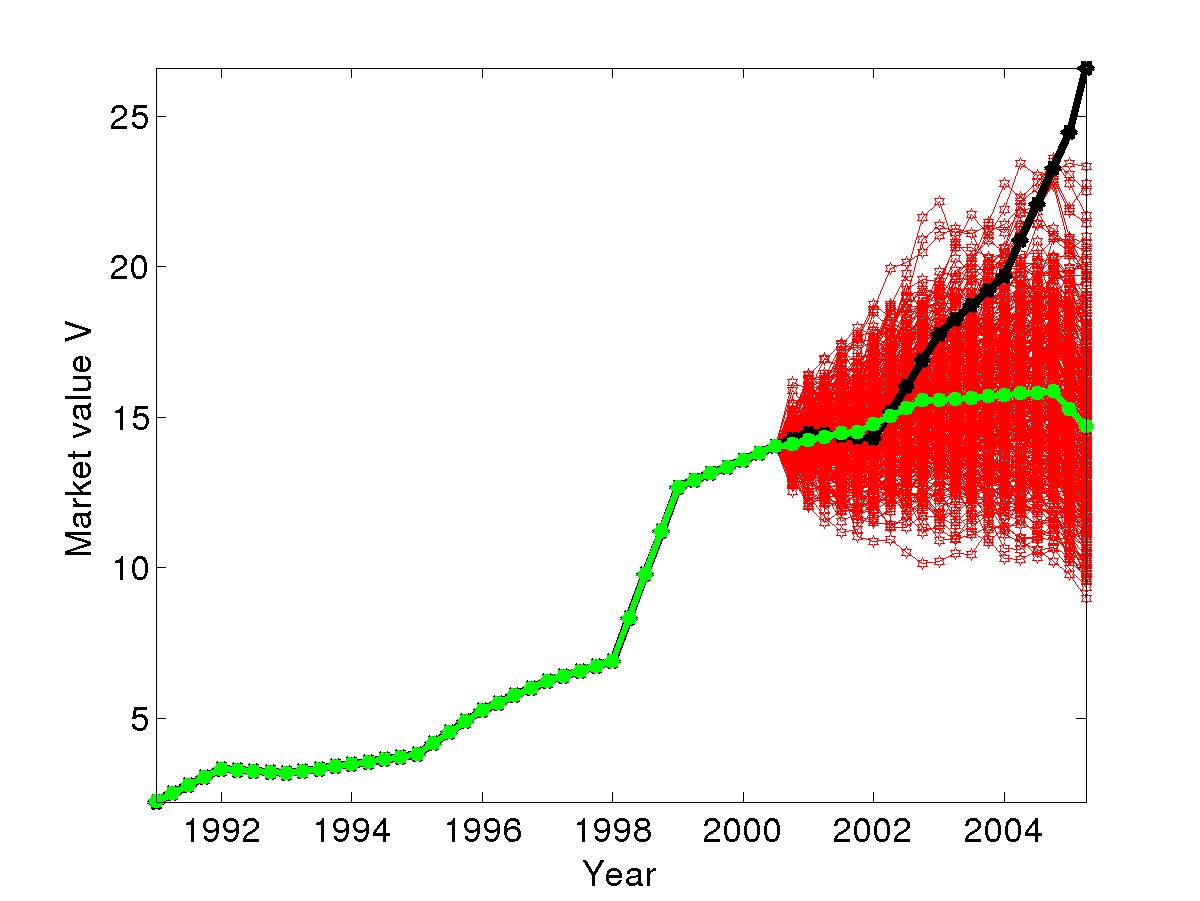}}
   \hskip 0.01\textwidth
   \subfigure[]{
   \label{FIG02d}
   \includegraphics[width=0.45\textwidth]{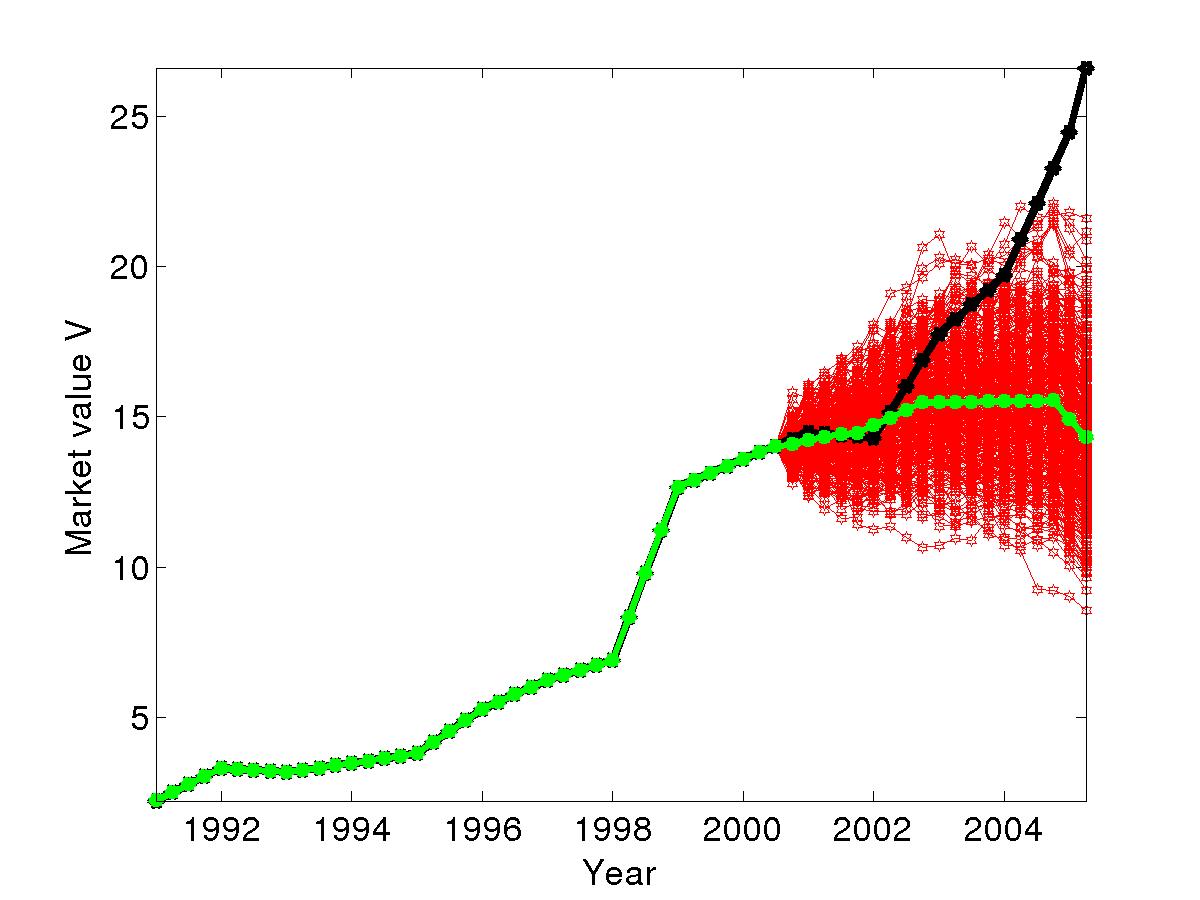}}
   \caption{The graphs at the left ((a) and (c) respectively for corporates $C_1$ and $C_2$) correspond to the delay model  while the graphs at the right ((b) and (d) respectively for corporates $C_1$ and $C_2$) 
correspond to  Merton model. We take $T=5\;\;L=9.5$ and the function $g$ is the quadratic interpolation of the 
 standard deviation of daily returns $\sigma$ in the memory part. 
We have plotted 400 samples of the numerical solution along with the expectation (the means)
 of the numerical solution (green curves). The curves of the real data of the firm value  $V$ as a function of time  are in black (black thick curves).}
 \label{FIG02}
\end{figure}

The following firms are used:
\begin{itemize}
 \item [$C_1$] {\tt Great Northern Iron Ore Pptys} (\figref{FIG02a} and \figref{FIG02b})
 \item [$C_2$] {\tt Tor Minerals Intl Inc }  (\figref{FIG02c} and \figref{FIG02d})
\item [$C_3$] {\tt South Jersey Inds Inc } (\figref{FIG011a} and \figref{FIG011b})
\item [$C_4$] {\tt Rentech  Inc} (\figref{FIG01c} and \figref{FIG01d})
\item [$C_5$] {\tt Magna International Inc} (\figref{FIG01a} and \figref{FIG01b})
\item [$C_6$]  {\tt First Citizens Bancshares Inc NC} (\figref{FIG011c} and \figref{FIG011d})
\end{itemize}


   As we have already mentioned, the time origin corresponds to the year (2000+1/2), the data before
are memory data and we want to predict the data after (2000+1/2).
We  plot 400 samples of the numerical solution for our delayed model and Merton model along with the means
 of the numerical samples (green curves). As the origin is year (2000+1/2), the part of the mean curves
 before the origin are just the curves of the real firm market value $V$  in that interval.
The curves of the  real firm market value $V$ as a function of time  are in black (black thick curves).
Indeed we want the means of the numerical samples (green curves) to fit well the  real firm market value $V$ (black thick curves) with moderate standard derivations 
(few spread of the numerical samples comparing to its mean), this will be the aim  of our comparisons.

In all  graphs, the function $g$ (volatility in delayed model) is the quadratic interpolation 
of the  standard deviation of daily returns $\sigma$ in the memory part while the volatility in the Merton model is just the mean of the memory part.

In \figref{FIG01} we take $L=T=9.5$, the graphs at the left hand  size.
(\figref{FIG01a} and \figref{FIG01c} respectively for firms $C_5$ and $C_4$ ) correspond to the delayed model while
the graphs at the right hand size (\figref{FIG01b} and \figref{FIG01d} respectively for firms $C_5$ and $C_4$) correspond to Merton model. 
For  corporate $C_{5}$, \figref{FIG01a} shows the good prediction with reasonable standard deviation (as the numerical samples are not much spread) of the delay model while \figref{FIG01b}  
shows the early good prediction of the Merton model but the prediction has failled just after the year 2005. For corporate $C_{4}$, \figref{FIG01b} and \figref{FIG01d} show the 
good prediction  before  2005 with relatively large standard deviation of the delayed model  and Merton model.

In \figref{FIG011}, we  aslo take $L=T=9.5$, the graphs at the left hand  size 
(\figref{FIG011a} and \figref{FIG011c} respectively for firms $C_{3}$ and $C_6$) correspond to the delayed model while
the graphs at the right hand size (\figref{FIG011b} and \figref{FIG011d}  respectively for firms $C_{3}$ and $C_6$)) correspond to Merton model.
For  corporate $C_{3}$, \figref{FIG011a} shows an early accepted prediction for delayed model comparing to the Merton model in \figref{FIG011b} where the prediction 
is more bad (black thick curve been really far away from green curve). For corporate $C_{6}$, we have the same observation as for corporate $C_{3}$ according to  \figref{FIG011c} and \figref{FIG011d}.

In \figref{FIG02} we take $L=9.5,\;T=5$, the graphs at the left (\figref{FIG02a} and \figref{FIG02c} respectively for firms $C_{1}$ and $C_2$) correspond to the delayed model while
the graphs at the right (\figref{FIG02b} and \figref{FIG02d} respectively for firms $C_{1}$ and $C_2$) correspond to Merton model. 
From \figref{FIG02}, we can observe that  for corporate $C_{1}$ the delayed model fit well the real data of the firm market value compared to the Merton model, while for corporate $C_2$ 
the two models fit well the real data of the firm market value before 2003.
\

\subsection{ Numerical Evaluation of Debt or Equity in a Levered Firm}
\subsubsection{Numerical schemes based on exponential integrators}
Debt and equity are linked by relation \eqref{en}, so we only need to solve one of the systems \eqref{eqou} and \eqref{eqouM}.
We consider here the random partial differential \eqref{eqou}, but where $C$,  $C_{y}$ and $r$ are time depending functions.
In our simulation we consider those values as   
time depending functions, which are constant within each year interval as we have in our data set. 
Indeed to solve  numerically this equation the domain of $v$ need to be troncated.  
Taking in  to account the fact that $C$,  $C_{y}$ and $r$ are time depending functions, we therefore have
\begin{eqnarray}
\label{eq8n}
\left\lbrace \begin{array}{l}
        \dfrac{1}{2}\,g^2(V(t-L))v^2\,f_{vv}+(r(t)v-C(t))f_v+f_t-r(t)f-C_y(t)+C(t)=0,  \\
f(v,T)=\max(v-B,0), \;\;\;\;\; v\in [0, V_{\max}]\\
f(0,t)=0, \;\quad \quad \;\;\;t \in [T-L,T]\\
f(V_{\max},t) =V_{\max}-B e^{-\int_{t}^{T}r(s)ds},\;\quad \quad \;\;\;t \in [T-L,T]
 \end{array} \right.
\end{eqnarray}
Our model problem \eqref{eq8n} is similar to Europeans call options prices, we can therefore  take $V_{\max}$ three or four times $B$ according to \cite{optionbound}.
In our simulation we take $V_{\max}=4B$, where $B$ is the amount that the firm must pay to the debtholders at the maturity date $T$
(like the strike price for options prices). Our system \eqref{eq8n} is a  backward  system, to transform it to the forward one,
  we use the tranformation $\tau=T-t$, and the corresponding equation is given by 
\begin{eqnarray}
\label{eq8nn}
\left\lbrace \begin{array}{l}
        \dfrac{1}{2}\,g^2(V(T-\tau-L))v^2\,f_{vv}+(r(\tau)v-C(\tau))f_v-r(\tau)f\\
       \qquad \qquad \qquad \qquad \qquad \qquad \qquad\qquad \qquad  -C_y(\tau)+C(\tau)=f_\tau ,  \\
f(v,0)=\max(v-B,0), \;\;\;\;\; v\in [0, V_{\max}]\\
f(0,\tau)=0,\;\quad \quad \;\;\tau \in [0,L] \\
f(V_{\max},\tau) =V_{\max}-B e^{-\int_{T-\tau}^{T}r(s)ds},\;\quad \quad \;\;\;\tau\in [0,L]
 \end{array} \right.
\end{eqnarray}
Please note that after the transformation  $\tau=T-t$,  the functions $C(\tau), C_{y}(\tau)$ and $r(\tau)$ in \eqref{eq8nn}  are normally the functions $C(T-\tau), C_{y}(T-\tau)$ and $r(T-\tau)$. 

To apply sophistical technique  to the convection term (the term with $f_{v}$) in order to avoid numerical instabilities, 
let us put this term in the so called the conservation form. 
In fact 
$$  (r(\tau)v-C(\tau))f_v= \left(r(\tau)v-C(\tau))f\right)_{v}- r(\tau) f$$
Using this relation, equation \eqref{eq8nn} become
\begin{eqnarray}
\label{eq8nnn}
\left\lbrace \begin{array}{l}
        \dfrac{1}{2}\,g^2(V(T-\tau-L))v^2\,f_{vv}+\left((r(\tau)v-C(\tau))f\right)_{v}-2r(\tau)f,  \\
        \qquad \qquad \qquad \qquad \qquad \qquad \qquad\qquad \qquad  -C_y(\tau)+C(\tau)=f_\tau\\
f(v,0)=\max(v-B,0), \;\;\;\;\; v\in [0, V_{max}]\\
f(0,\tau)=0,\;\quad \quad \;\;\tau \in [0,L] \\
f(V_{max},\tau) =V_{max}-B e^{-\int_{T-\tau}^{T}r(s)ds},\;\quad \quad \;\;\;\tau\in [0,L]
 \end{array} \right.
\end{eqnarray}
To solve equation \eqref{eq8nnn} two cases can be considered:
\begin{enumerate}
 \item  The case where $T-\tau-L \leq 0,\;\;\; \forall \tau \in [0,L]$, then $T\leq L$.
\item The case where  $T>L$.
\end{enumerate}
For the first case ($T \leq L$) the RPDE \eqref{eq8nnn} become the deterministic PDE
since $V(t) =\varphi(t)\;\;$ for $ t\in [-L,0]$ as given in \eqref{model}.

For  the second case ($T>L$), to solve  \eqref{eq8nnn} the following step should be followed
 \begin{enumerate}
  \item Solve the stochastic equation \eqref{model} to have a sample of the numerical solution of $V$ as we did in  the previous section.\\
   \item  Use the numerical sample solution of $V$ from step 1 to build the diffusion coefficient (the coefficient of $f_{vv}$) in the RPDE \eqref{eq8nnn},
 which therefore become a deterministic PDE for this fixed numerical sample of $V$.
 \item Solve  the deterministic PDE from step 2 for the fixed numerical sample of $V$ from step 1.
\item Repeat step 1, step 2 and step 3, $M$ times (relatively large) and use the Monte Carlo technique to estimate the expectation value of $f$ 
and also any moment of the stochastic process $f$ if need.
 \end{enumerate}
As the two cases require the solution of the deterministic PDE, in the sequel we  will consider the first case ($T\leq L$), and  the corresponding
deterministic PDE is given by
\begin{eqnarray}
\label{eq8nnnn}
\left\lbrace \begin{array}{l}
        \dfrac{1}{2}\,g^2(\varphi(T-\tau-L))v^2\,f_{vv}+\left((r(\tau)v-C(\tau))f\right)_{v}-2r(\tau)f\\
        \qquad \qquad \qquad \qquad \qquad \qquad \qquad\qquad \qquad -C_y(\tau)+C(\tau)=f_\tau ,  \\
f(v,0)=\max(v-B,0), \;\;\;\;\; v\in [0, V_{max}]\\
f(0,\tau)=0,\;\quad \quad \;\;\tau \in [0,L] \\
f(V_{max},\tau) =V_{max}-B e^{-\int_{T-\tau}^{T}r(s)ds},\;\quad \quad \;\;\;\tau\in [0,L].
 \end{array} \right.
\end{eqnarray}
For the discretization in the direction of $v$, we use  the combined finite difference--finite volume method. The interval $[0,V_{max}]$ is subdivised into $N$ parts
that we assume equal without loss the generality. As in center finite volume method, we approximate $f$ at the center of each interval. 
The diffusion part of the equation is approximated using the finite difference while the convection term is approximated using the standard
upwinding usual used in porous media flow problems \cite{ATthesis,Antoine,FV,TLGspe1}.

Let 
$$v_{i}=(2i-1)h/2,\;\;\; h=\dfrac{V_{max}}{N},\;\;\ i= 1,2,.....,N$$ being the center of each subdivision .
We approximate the diffusion  term at each center by 
\begin{eqnarray*}
 \dfrac{1}{2}\,g^2(\varphi(T-\tau-L))v_{i}^2\,f_{vv}(v_{i})&\approx&\dfrac{1}{2 h^{2}}\,g^2(\varphi(T-\tau-L))v_{i}^2 \times\\  && \left( f_{i+1}(\tau)-2f_{i}(\tau)+ f_{i-1}(\tau)\right),\;i = 2,..., N-1.
\\
\dfrac{1}{2}\,g^2(\varphi(T-\tau-L))v_{1}^2\,f_{vv}(v_{1})&\approx& \dfrac{2}{3h}\,g^2(\varphi(T-\tau-L))v_{1}^2 \times\\  && \left( \dfrac{f_{2}(\tau)-f_{1}(\tau)}{h}-2\dfrac{f_{1}(\tau)}{h}\right)\\
\dfrac{1}{2}\,g^2(\varphi(T-\tau-L))v_{N}^2\,f_{vv}(v_{N})&\approx& \dfrac{2}{3h}\,g^2(\varphi(T-\tau-L))v_{N}^2 \times\\
 &&  \left( \dfrac{f(V_{max},\tau)-f_{N}(\tau)}{h/2}-\dfrac{f_{N}(\tau)-f_{N-1}(\tau)}{h}\right)
\end{eqnarray*}

 This approximation is similar to the one in \cite{option}  with central difference on  non uniform grid.
We approximate the convection term using the standard upwinding technique as following
\begin{eqnarray*}
 \left(r(\tau)v-C(\tau))f\right)_{v}(v_{i})\approx \dfrac{(r(\tau)v_{i+1/2}-C(\tau))f_{i}^{+}(\tau)-(r(\tau)v_{i-1/2}-C(\tau))f_{i-1}^{+}(\tau)}{h}
\end{eqnarray*}
where
\begin{eqnarray}
 f_{i}^{+}(\tau)&=& \left\lbrace \begin{array}{l}
   f_{i} (\tau)\;\;\;\;\; \text{if}\;\;\;\;r(\tau)v_{i+1/2}-C(\tau)\geqslant 0\;\;\;\;\;\\
     f_{i-1}(\tau)  \;\;\;\;\; \text{if}\;\;\;\;r(\tau)v_{i+1/2}-C(\tau)<0\;\;\;\;\;      
             \end{array}\right.\\
    v_{i+1/2}&=&v_{i}+h/2,\;\;\;\; v_{i-1/2}=v_{i}-h/2=v_{i-1}+h/2,
\end{eqnarray}
where 
$$ f_{i}(\tau) \approx f(v_{i},\tau).$$

Reorganizing  all  previous diffusion and convection approximations  lead to the  following initial value problem 
\begin{eqnarray}
\label{ode}
\left\lbrace \begin{array}{l}
 \dfrac{d\mathbf{f}}{d\tau}=\mathbf{A}(\tau)\mathbf{f}+\mathbf{b}(\tau),\;\;\;\;\;\; \tau\in [0,L]\\
\mathbf{f}(0)=\left(\max(v_{1}-B,0),...,\max(v_{N}-B,0)\right)^{T}.
\end{array}\right.
\end{eqnarray}
where $\mathbf{A}(\tau)$ is a tridiagonal matrix  and
\begin{eqnarray}
 \mathbf{f}(\tau)=\left( f_{i}(\tau)\right)_{1\leq i\leq N}, \mathbf{b}(\tau) = C(\tau)-C_{y}(\tau) + \mathbf{k}(\tau).
\end{eqnarray}
 where $\mathbf{k}$ is the contribution from boundary conditions.

The function $x\mapsto \max(x,0)$ is not smooth, it important to approximate it by a smooth function.
The approximation in \cite{option} is a fourth-order smooth function denoted $\pi_{\epsilon}$ and defined by 
\begin{eqnarray}
 \pi_{\epsilon}(x)=\left\lbrace \begin{array}{l}
                              x \;\;\;\; \text{if}\;\;\;\;\;\;x\geqslant \epsilon\\
                               c_{0}+c_{1}x+.....+c_{9}x^{9}   \;\;\;\; \text{if}\;\;\;\;\;\; -\epsilon< x<\epsilon\\
                               0  \;\;\;\; \text{if}\;\;\;\;\;  x \leq -\epsilon
                             \end{array}\right.
\end{eqnarray}
where  $0<\epsilon \ll 1$ is the transition parameter and
\begin{eqnarray*}
c_{0}=\dfrac{35}{256}\epsilon,\;\;\;\;c_{1}=\dfrac{1}{2},\;\;\;c_{2}=\dfrac{35}{64 \epsilon},\;\;\;\;c_{4}=-\dfrac{35}{128\epsilon^{3}}
,\,\,\,\\
c_{6}=\dfrac{7}{64\epsilon^{5}},\,\,c_{8}=-\dfrac{5}{256\epsilon^{7}},\;\; c_{3}=c_{5}=c_{7}=c_{9}=0. 
\end{eqnarray*}
This approximation allow us to write
\begin{eqnarray}
 \mathbf{f}(0)= \pi_{\epsilon}(\mathbf{v}-B),\;\;\; \mathbf{v}= \left( v_{i}\right)_{1\leq i\leq N}.
\end{eqnarray}
 Let us introduce the time stepping discretization for the ODE \eqref{ode} based on exponential integrators. 
Classical numerical methods usually used  are Implicit Euler scheme and Crank--Nicolson scheme \cite{P3}.
Following works from \cite{option,ATthesis,Antoine} the exact solution of \eqref{ode} is given by
\begin{eqnarray}
\label{mild}
 \mathbf{f}(\tau_{n}+\Delta \tau)&=&e^{\int_{\tau_{n}}^{\tau_{n}+\Delta \tau}\mathbf{A}(s)ds} \left[\mathbf{f}(\tau_{n})+ \int_{\tau_{n}}^{\tau_{n}+\Delta \tau}e^{-\int_{\tau_{n}}^{s }\mathbf{A}(y)dy} \mathbf{b}(s)ds \right]
\\ &&\qquad \qquad \tau_{n}= n\,\Delta \tau,\;\;\;\;\; n=0,..., M,\;\;\;\;\ \Delta \tau >0. 
\end{eqnarray}

Note that \eqref{mild} is the exact representation of the solution, to have the numerical schemes, approximations are needed, the first approximations (using the quadrature rule) may be 
\begin{eqnarray}
\label{app}
 \int_{t_{n}}^{t_{n}+\Delta \tau}\mathbf{A}(s)ds\approx \Delta \tau\mathbf{A}(\tau_{n})\;\;\;\;\int_{\tau_{n}}^{s }\mathbf{A}(y)dy\approx (s-\tau_{n})\mathbf{A}(\tau_{n})
\end{eqnarray}

Using these approximations we therefore have the following second-order approximation

\begin{eqnarray}
\label{milda}
 \mathbf{f}(\tau_{n}+\Delta \tau)&\approx &e^{\Delta \tau\mathbf{A}(\tau_{n})} \left[\mathbf{f}(\tau_{n})+ \int_{\tau_{n}}^{\tau_{n}+\Delta \tau}e^{-(s-\tau_{n})\mathbf{A}(\tau_{n})} \mathbf{b}(s)ds \right]
\end{eqnarray}
The simple scheme called  Exponential Differential scheme of order 1 (ETD1)
 is obtained by approximating  $\mathbf{b}(s)$  by  the constant $\mathbf{b}(\tau_n)$ and is given by 
\begin{eqnarray}
 \mathbf{f}_{n+1}=  \mathbf{f}_{n}+ \left(\Delta \tau\mathbf{A}(\tau_{n})\right)^{-1}\left[ \mathbf{A}(\tau_{n}) \mathbf{f}_{n} +\mathbf{b}(\tau_{n})
\right].
\end{eqnarray}
A second order scheme is given in \cite{option}.

Following the work in \cite[Lemma 4.1]{P4}, if the the function $\mathbf{b}$ can be well approximated by the polynomial of degree $p$
(which is the case here since we have  the exponential decay at the boundary $v=V_{max}$), from \eqref{milda}  we have   
\begin{eqnarray}
\label{hig}
 \mathbf{f}_{n+1}=\varphi_{0}(\Delta \tau\mathbf{A}(\tau_{n}))\mathbf{f}_{n} +\underset{j=0}{\sum^{p-1}}\underset{l=0}{\sum^{j}} 
\dfrac{\tau_{n}^{j-l}}{(j-l)!} \Delta \tau^{l+1} \varphi_{l+1} (\Delta \tau\mathbf{A}(\tau_{n}))\mathbf{b}_{j+1},
\end{eqnarray}
where
 \begin{eqnarray*}
  \mathbf{b}(\tau)&\approx&\underset{j=0}{\sum^{p-1}}\dfrac{\tau^{j}}{j!}\mathbf{b}_{j+1},\\
 \varphi_{0}(x)&=&e^{x},\;\;\;\varphi_{l}(x)=x\varphi_{l+1}(x)+\dfrac{1}{l!},\,\;\;\; l=0,1,2,.....
 \end{eqnarray*}
 Note that to have high  order  accuracy in time for $p>2$, the integral in \eqref{app} should be approximated more accurately.  The Magnus expansion may also used in such case.

 \subsubsection{Application with corporate data}
All schemes here can be implemented using Krylov subspace technique in the  computation the expomential 
functions presented in those schemes with the Matlab functions {\tt expmvp.m} or
{\tt phipm.m} from \cite{P4,P5}.  The Krylov subspace dimension we use is $m=10$ and the tolerance using in  the computation of the expomential
functions $\varphi_{i}$ is  $tol=1e-6$. We use $p=2$  and  obtain second order accuracy in time as the approximations \eqref{app} are second order in time. 

\begin{figure}[H]
 \subfigure[]{
   \label{FIG03a}
   \includegraphics[width=0.5\textwidth]{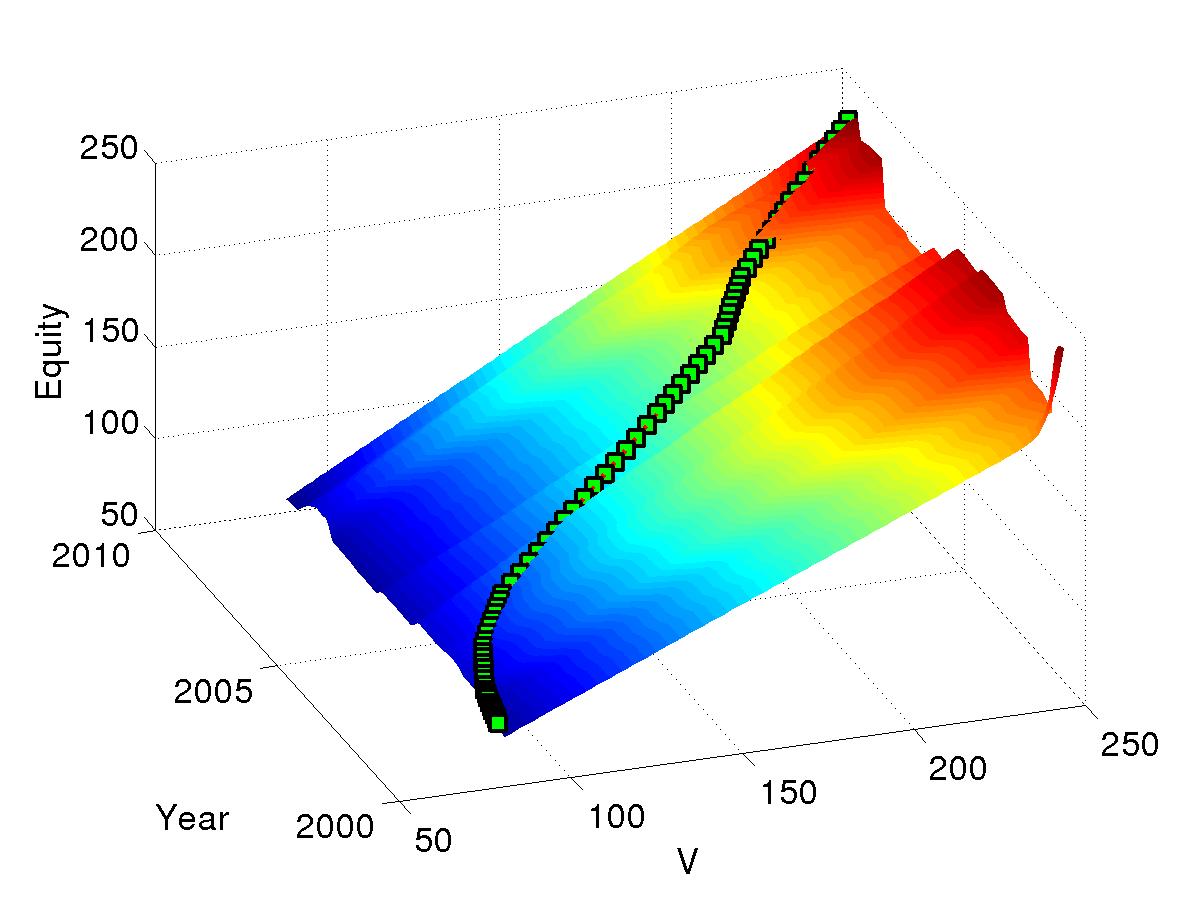}}
   \hskip 0.01\textwidth
   \subfigure[]{
   \label{FIG03b}
   \includegraphics[width=0.5\textwidth]{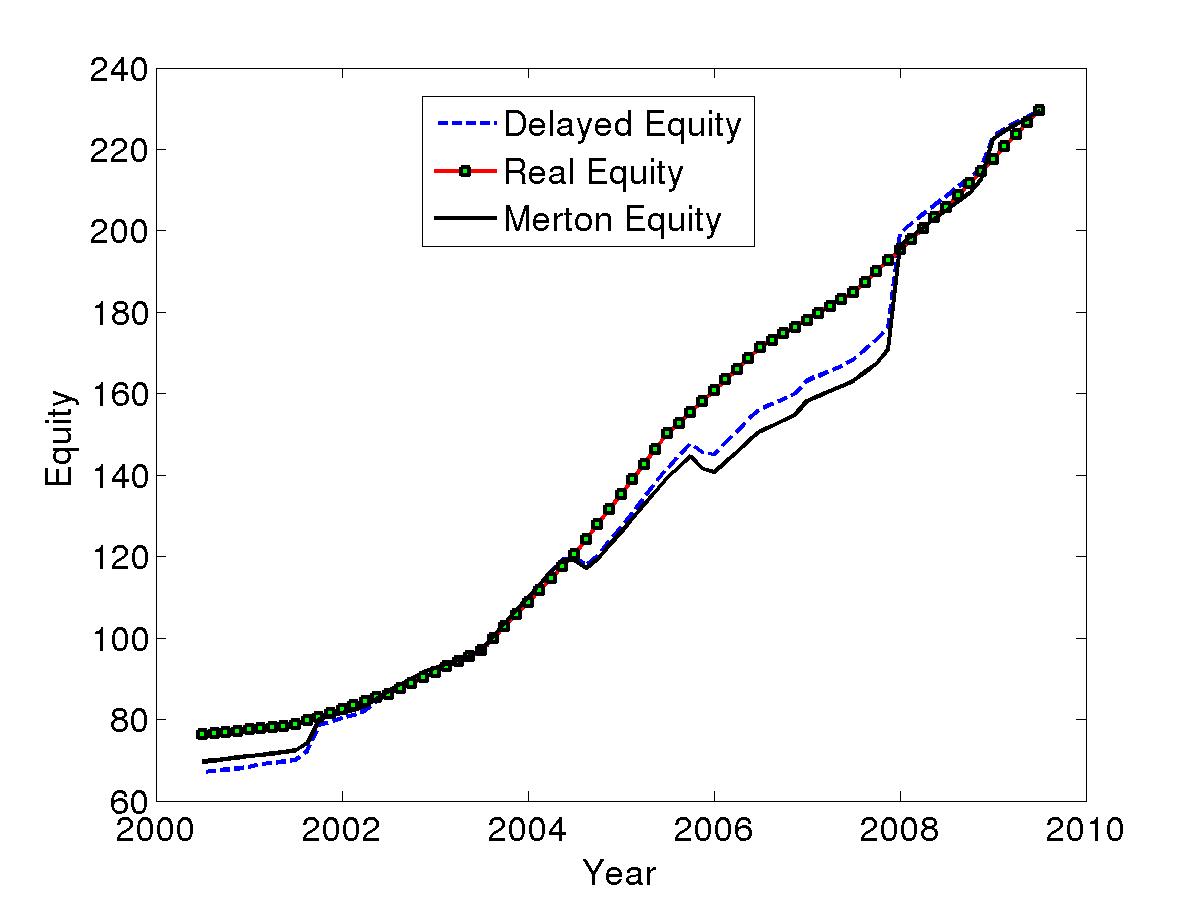}}
   \hskip 0.01\textwidth
   \subfigure[]{
   \label{FIG03c}
   \includegraphics[width=0.5\textwidth]{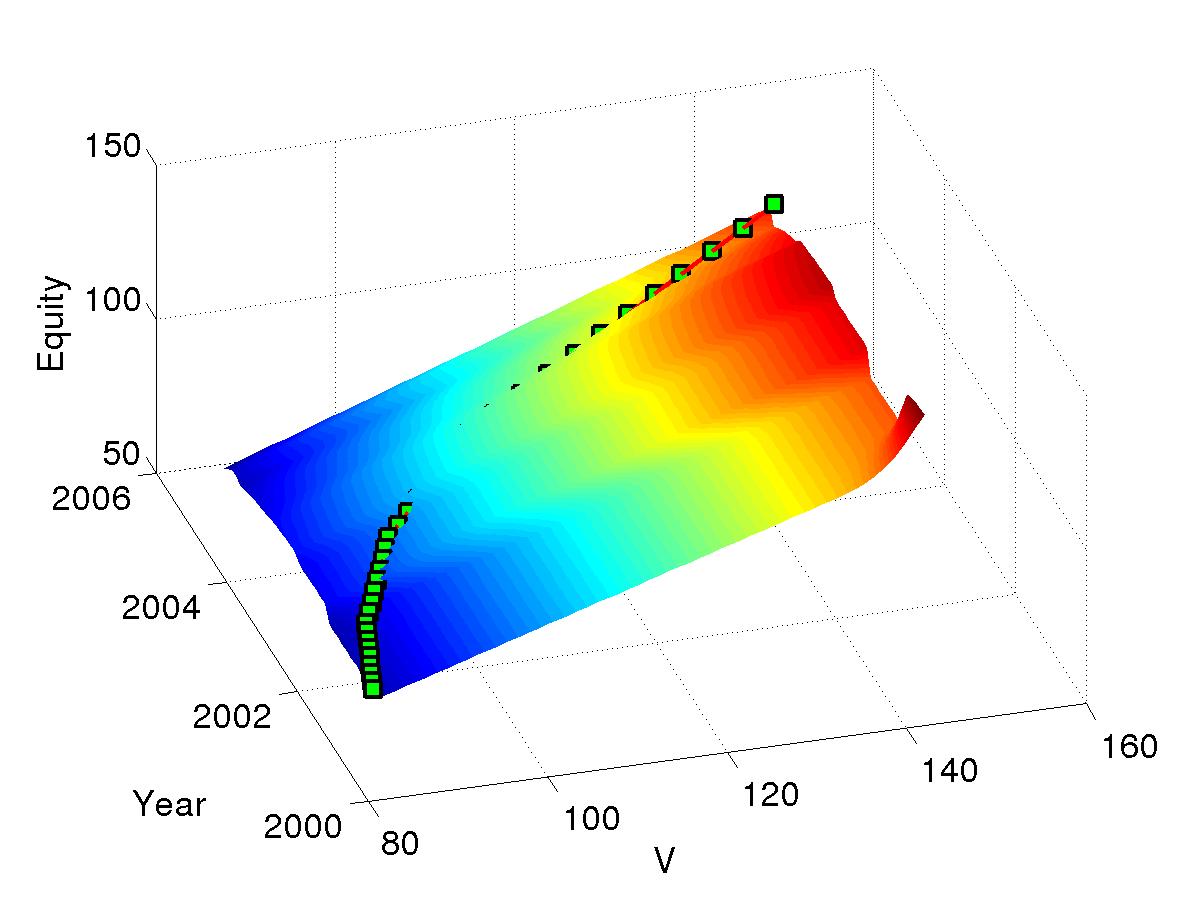}}
   \hskip 0.01\textwidth
   \subfigure[]{
   \label{FIG03d}
   \includegraphics[width=0.5\textwidth]{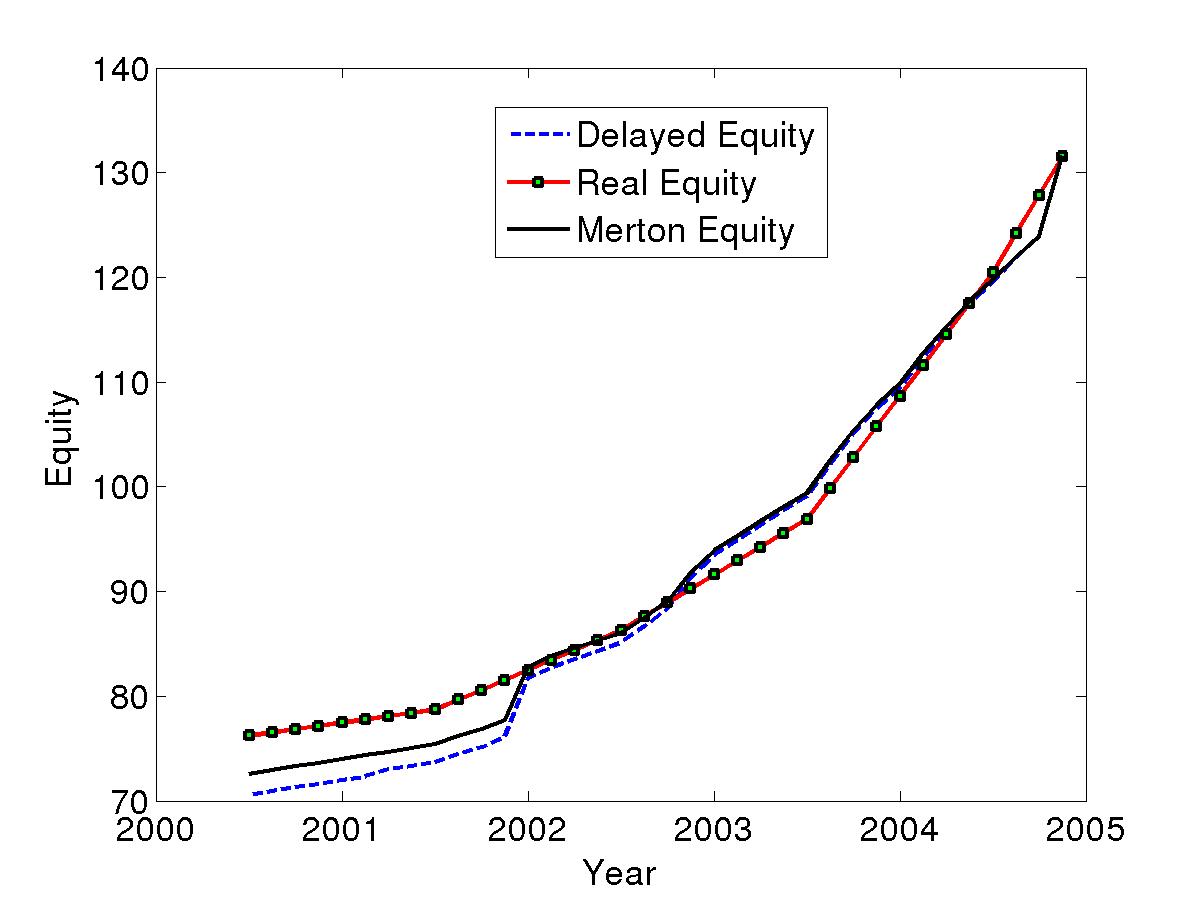}}
 \caption{ The graphs for  Firm $C_5$.  We plot at the left ((a) and (c) respectively)  the surface graphs of the numerical equity from our delayed model at $T=9.5$ and $T=5$.
 The corresponding  3 D graphs (green curves) of the real data of the firm equity value as a function of the time  and $V$ are also plotted in (a) and (c). 
 At the right ((b) and (d)), 2 D  graphs of the firm  equity  value as a function of time, corresponding to the surface graphs at  the left ((a) and (c)) respectively) are presented. Those 2 D graphs
 contain  the numerical equity from our delayed model, the numerical equity of the Merton model and the real data equity of the firm.}
 \label{FIG03}
\end{figure}

\begin{figure}[H]
 \subfigure[]{
   \label{FIG04a}
   \includegraphics[width=0.5\textwidth]{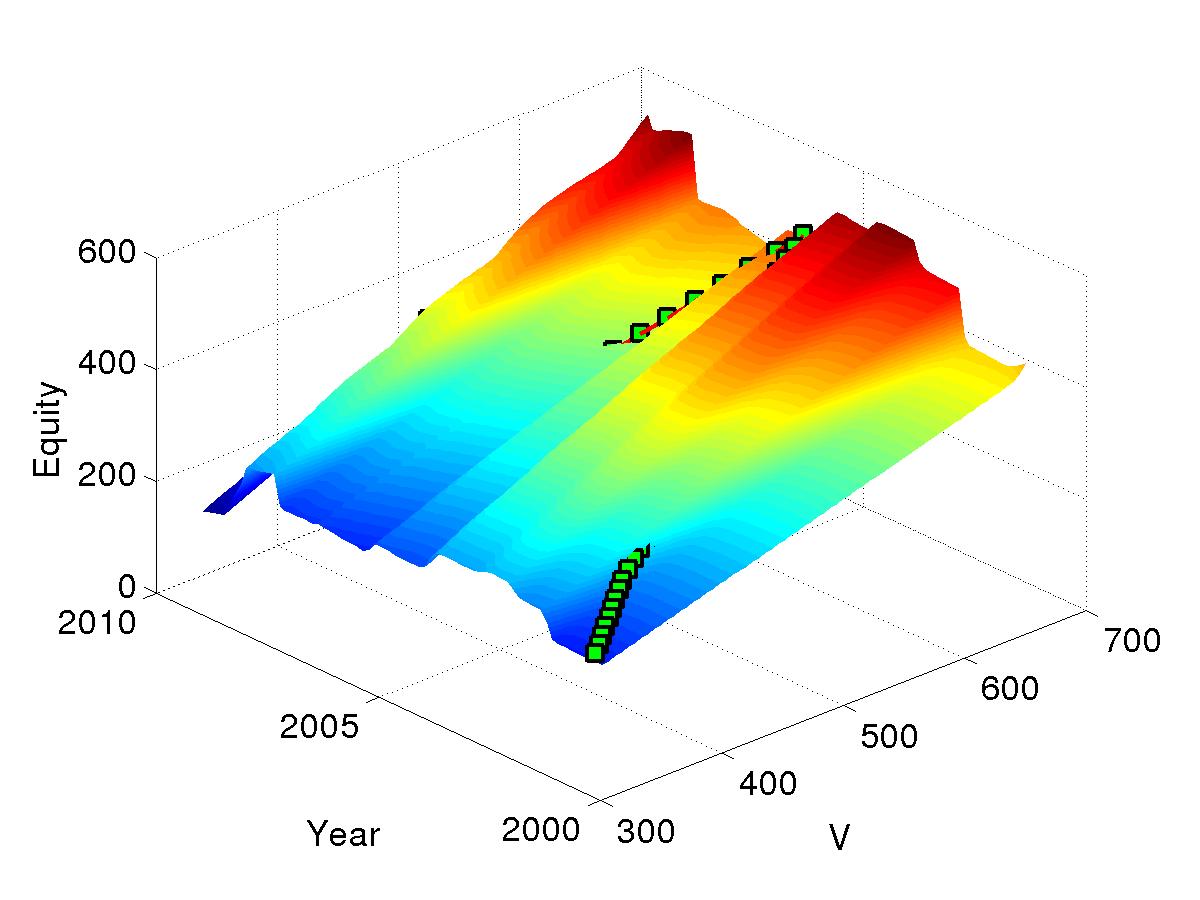}}
   \hskip 0.01\textwidth
   \subfigure[]{
   \label{FIG04b}
   \includegraphics[width=0.5\textwidth]{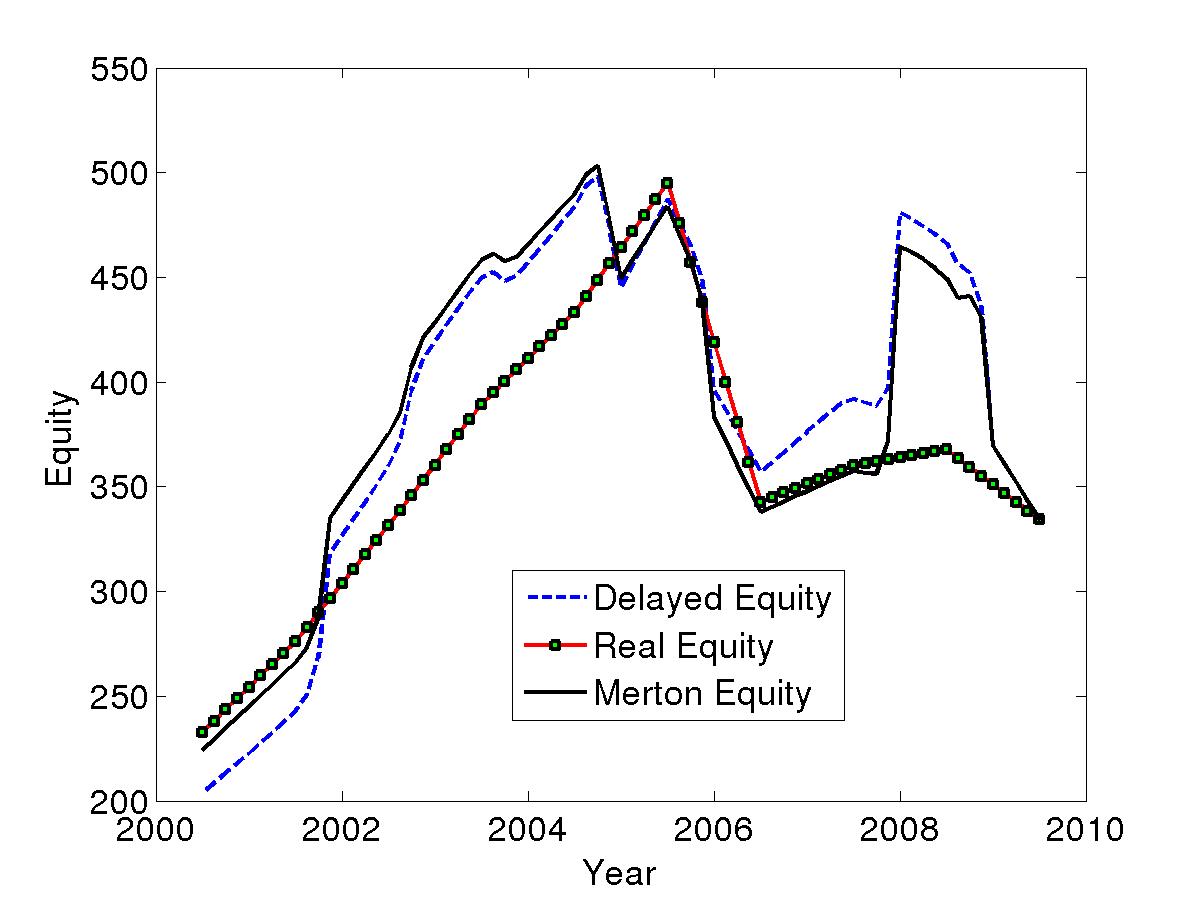}}
   \hskip 0.01\textwidth
   \subfigure[]{
   \label{FIG04c}
   \includegraphics[width=0.5\textwidth]{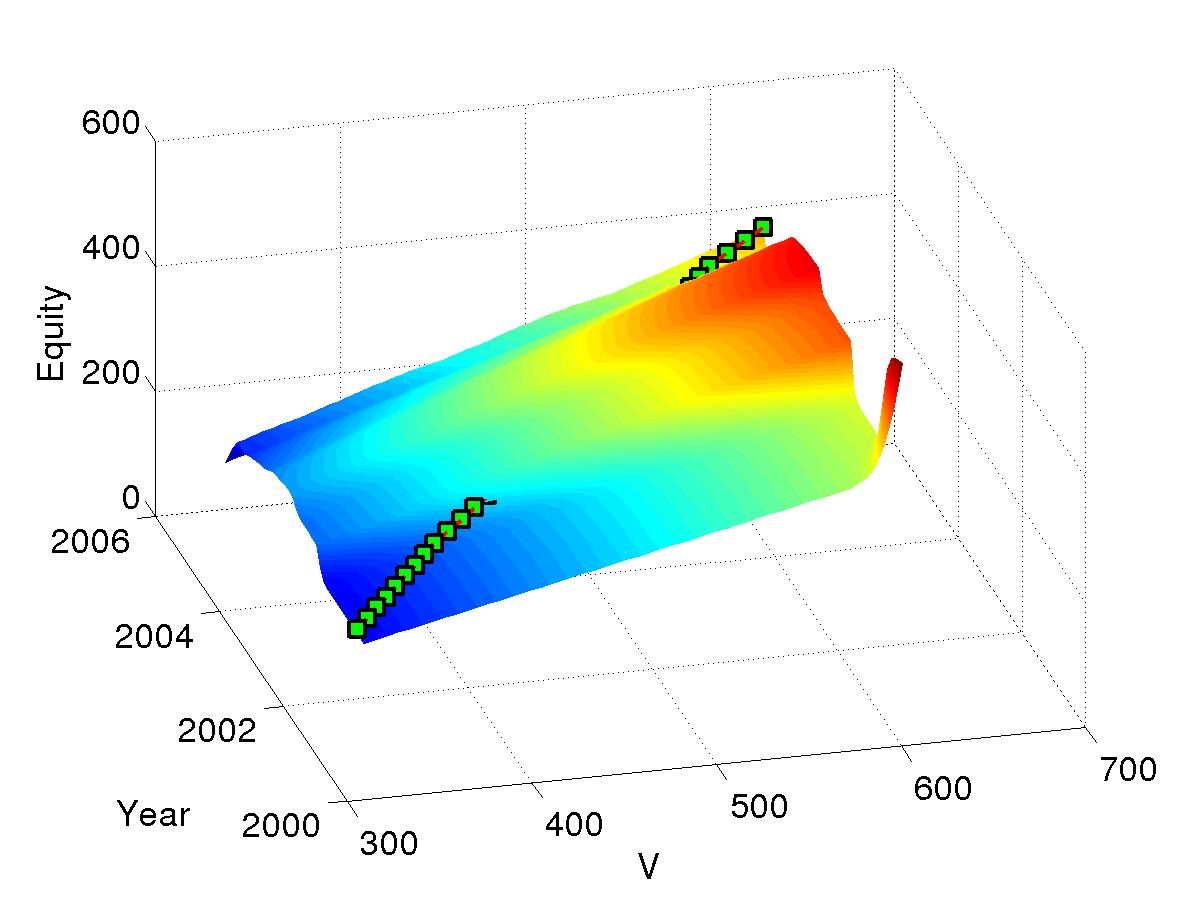}}
   \hskip 0.01\textwidth
   \subfigure[]{
   \label{FIG04d}
   \includegraphics[width=0.5\textwidth]{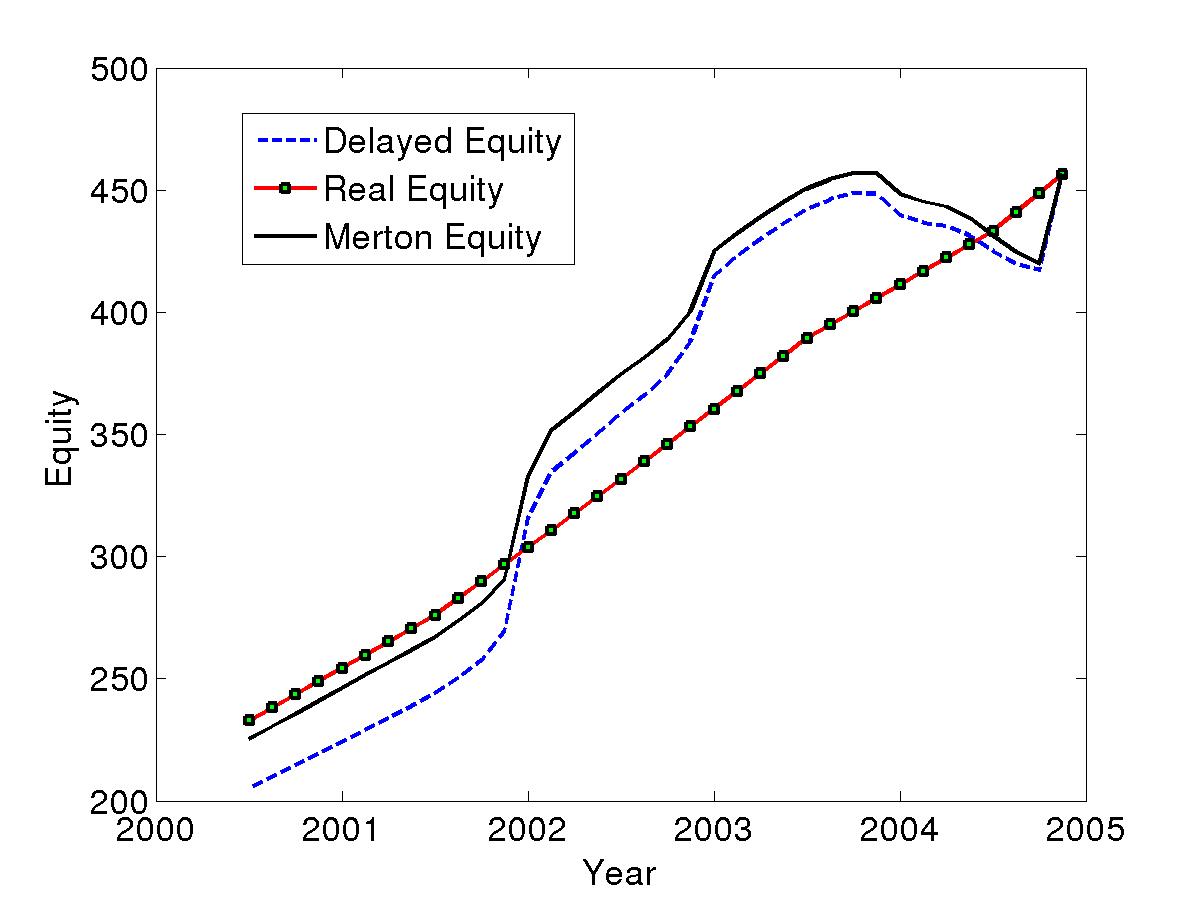}}
 \caption{ The graphs for  Firm $C_6$.  We plot at the left ((a) and (c) respectively)  the surface graphs of the numerical equity from our delayed model at $T=9.5$ and $T=5$.
 The corresponding  3 D graphs (green curves) of the real data of the firm equity value as a function of the time  and $V$ are also plotted in (a) and (c). 
 At the right ((b) and (d)), 2 D  graphs of the firm  equity  value as a function of time, corresponding to the surface graphs at  the left ((a) and (c)) respectively) are presented. Those 2 D graphs
 contain  the numerical equity from our delayed model, the numerical equity of the Merton model and the real data equity of the firm.}
 \label{FIG04}
\end{figure}

\begin{figure}[H]
 \subfigure[]{
   \label{FIG05a}
   \includegraphics[width=0.5\textwidth]{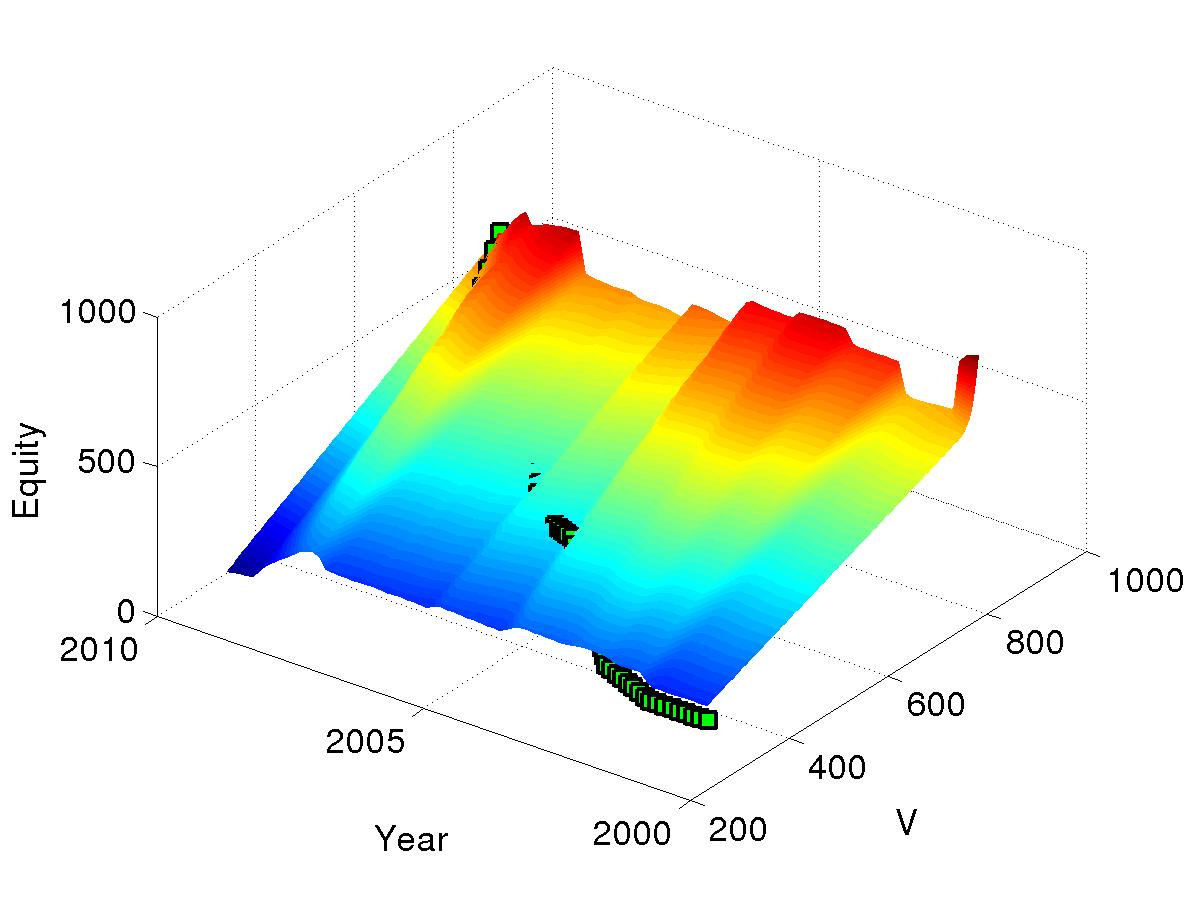}}
   \hskip 0.01\textwidth
   \subfigure[]{
   \label{FIG05b}
   \includegraphics[width=0.5\textwidth]{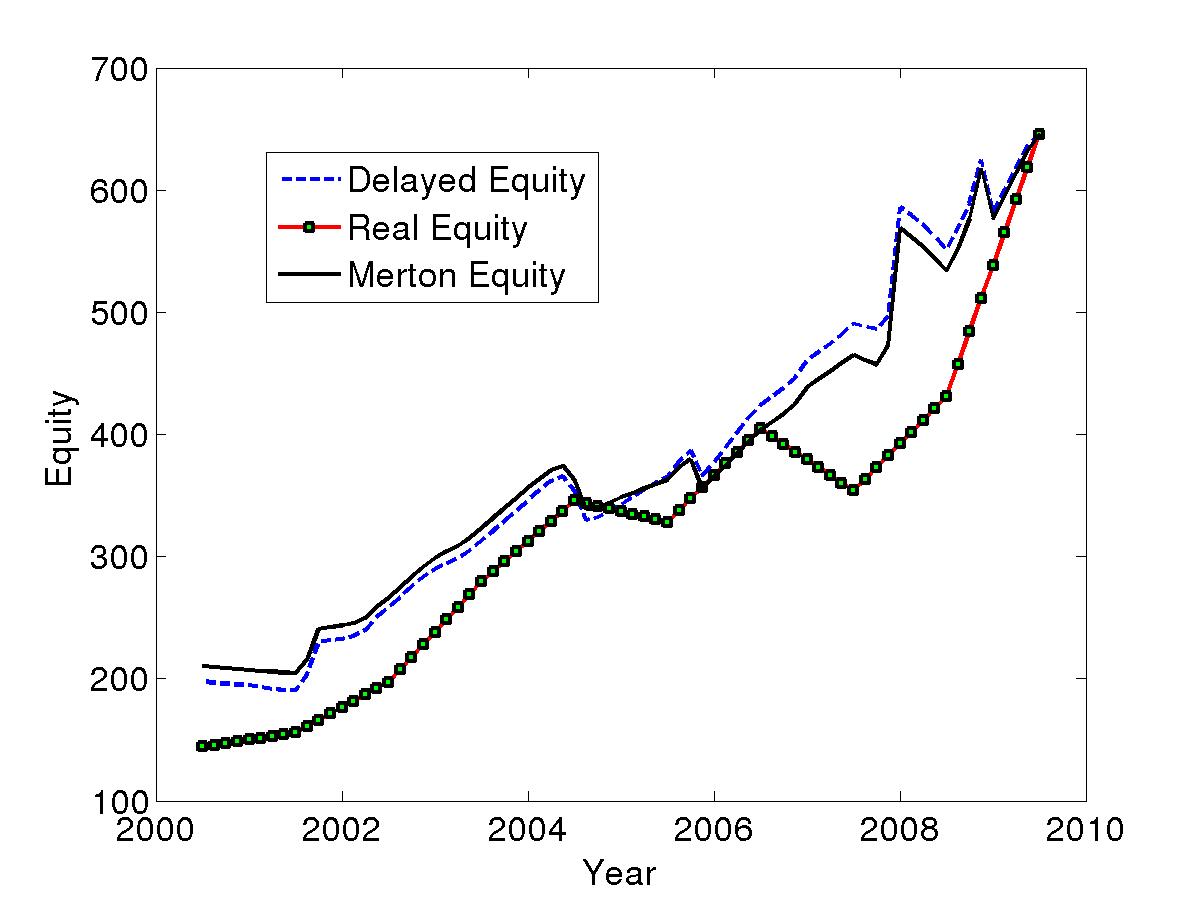}}
   \hskip 0.01\textwidth
   \subfigure[]{
   \label{FIG05c}
   \includegraphics[width=0.5\textwidth]{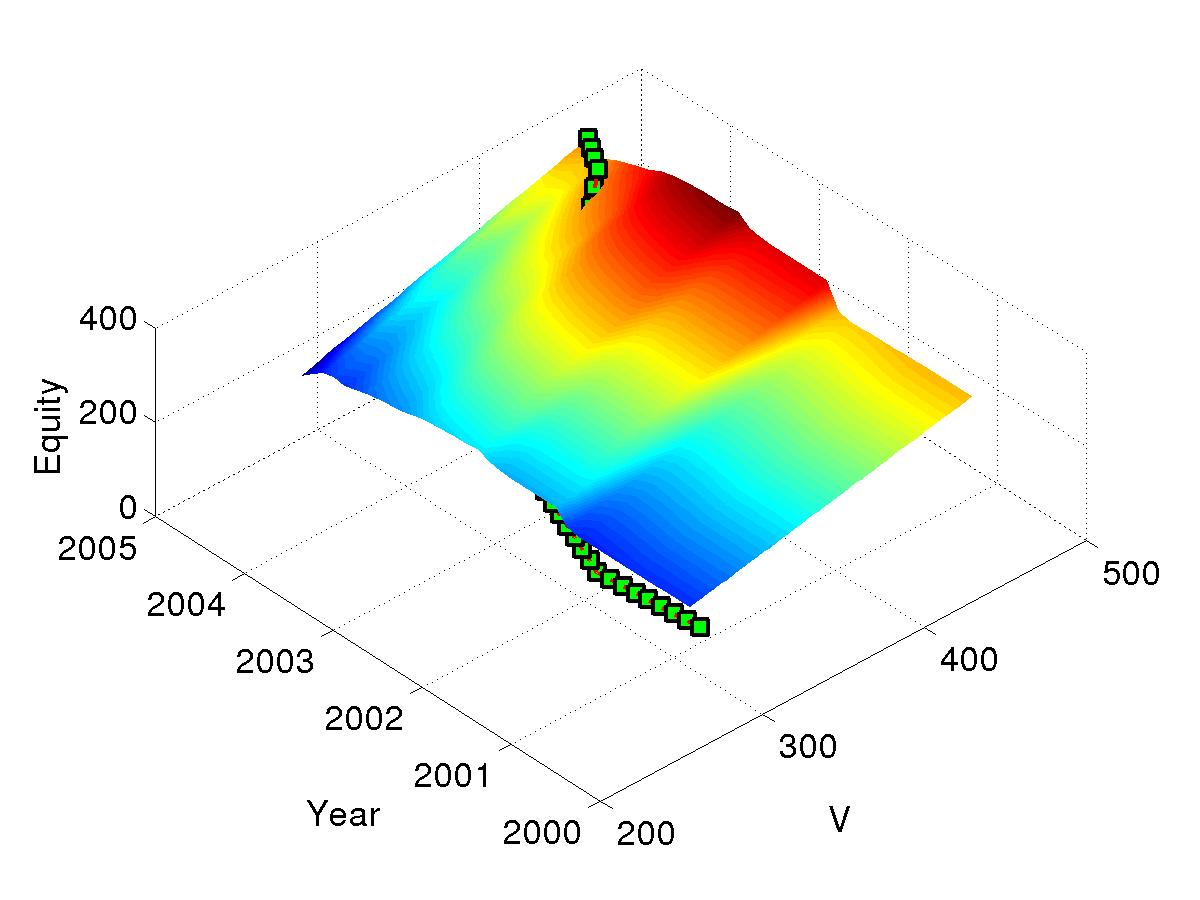}}
   \hskip 0.01\textwidth
   \subfigure[]{
   \label{FIG05d}
   \includegraphics[width=0.5\textwidth]{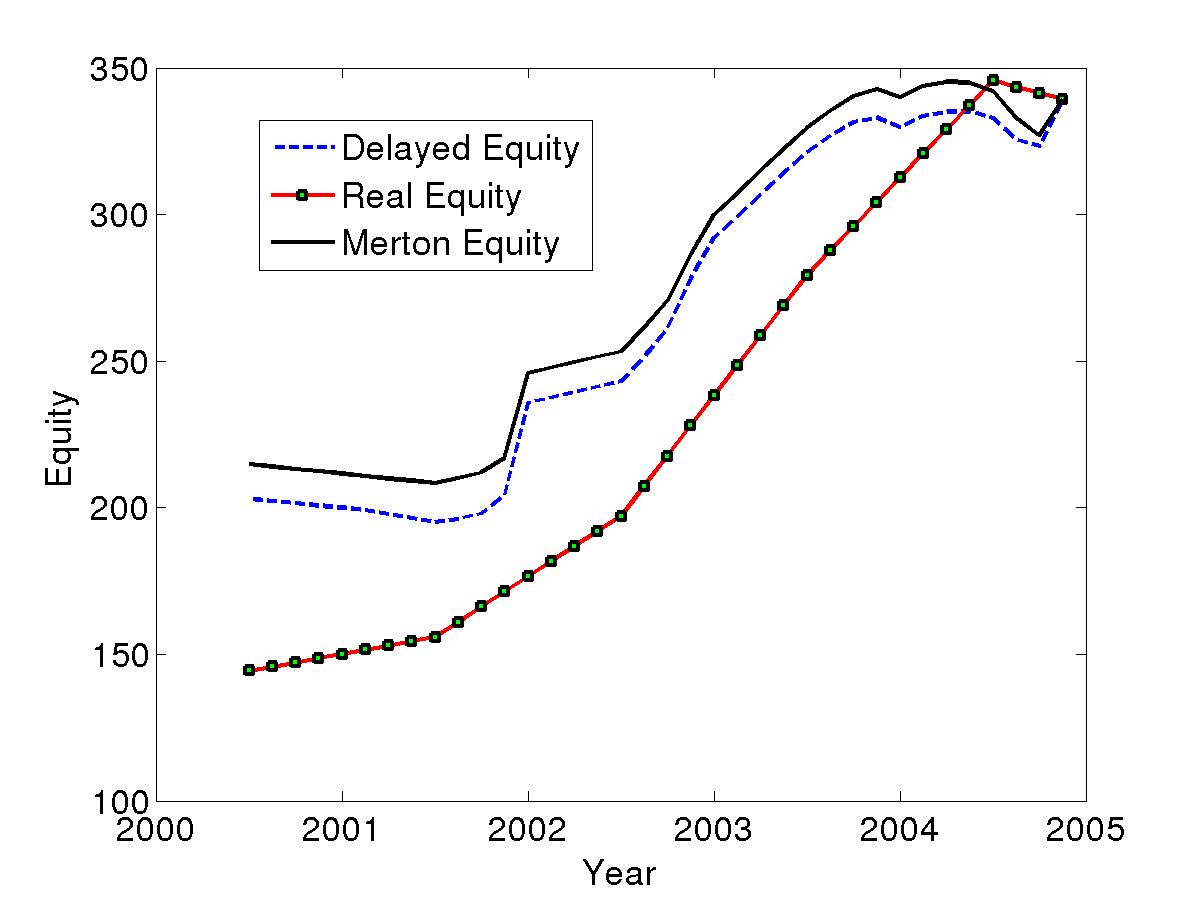}}
 \caption{ The graphs for  Firm $C_7$.  We plot at the left ((a) and (c) respectively)  the surface graphs of the numerical equity from our delayed model at $T=9.5$ and $T=5$.
 The corresponding  3 D graphs (green curves) of the real data of the firm equity value as a function of the time  and $V$ are also plotted in (a) and (c). 
 At the right ((b) and (d)), 2 D  graphs of the firm  equity  value as a function of time, corresponding to the surface graphs at  the left ((a) and (c)) respectively) are presented. Those 2 D graphs
 contain  the numerical equity from our delayed model, the numerical equity of the Merton model and the real data equity of the firm.}
 \label{FIG05}
\end{figure}

\begin{figure}[H]
 \subfigure[]{
   \label{FIG06a}
   \includegraphics[width=0.5\textwidth]{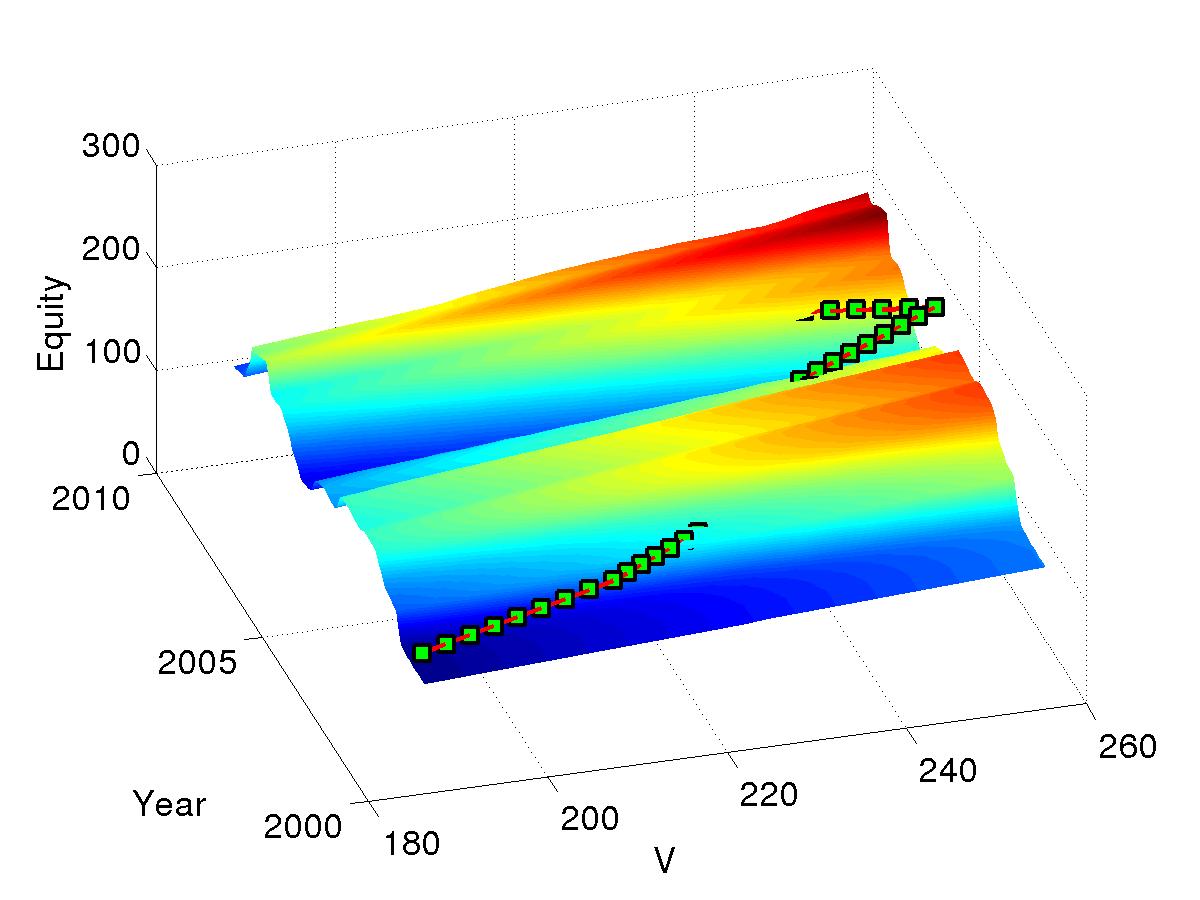}}
   \hskip 0.01\textwidth
   \subfigure[]{
   \label{FIG06b}
   \includegraphics[width=0.5\textwidth]{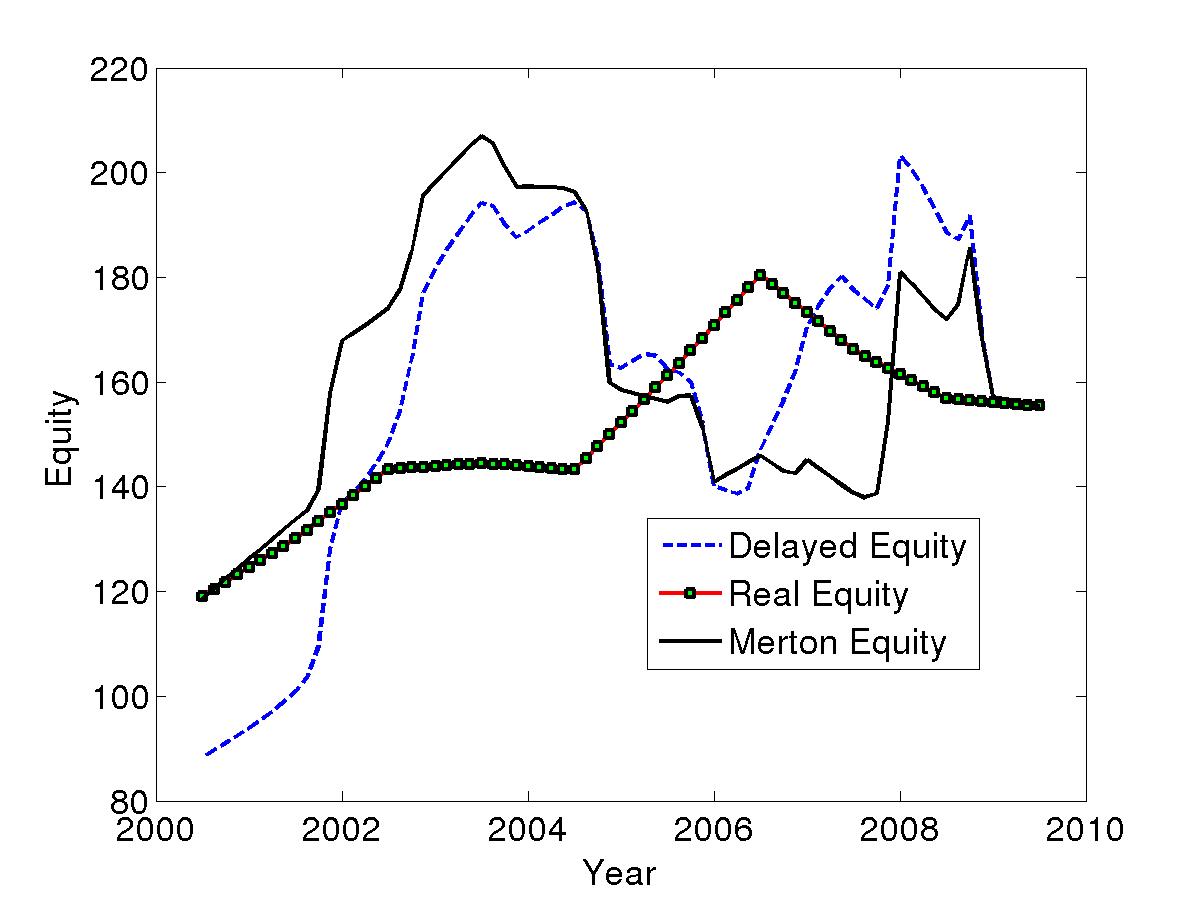}}
   \hskip 0.01\textwidth
   \subfigure[]{
   \label{FIG06c}
   \includegraphics[width=0.5\textwidth]{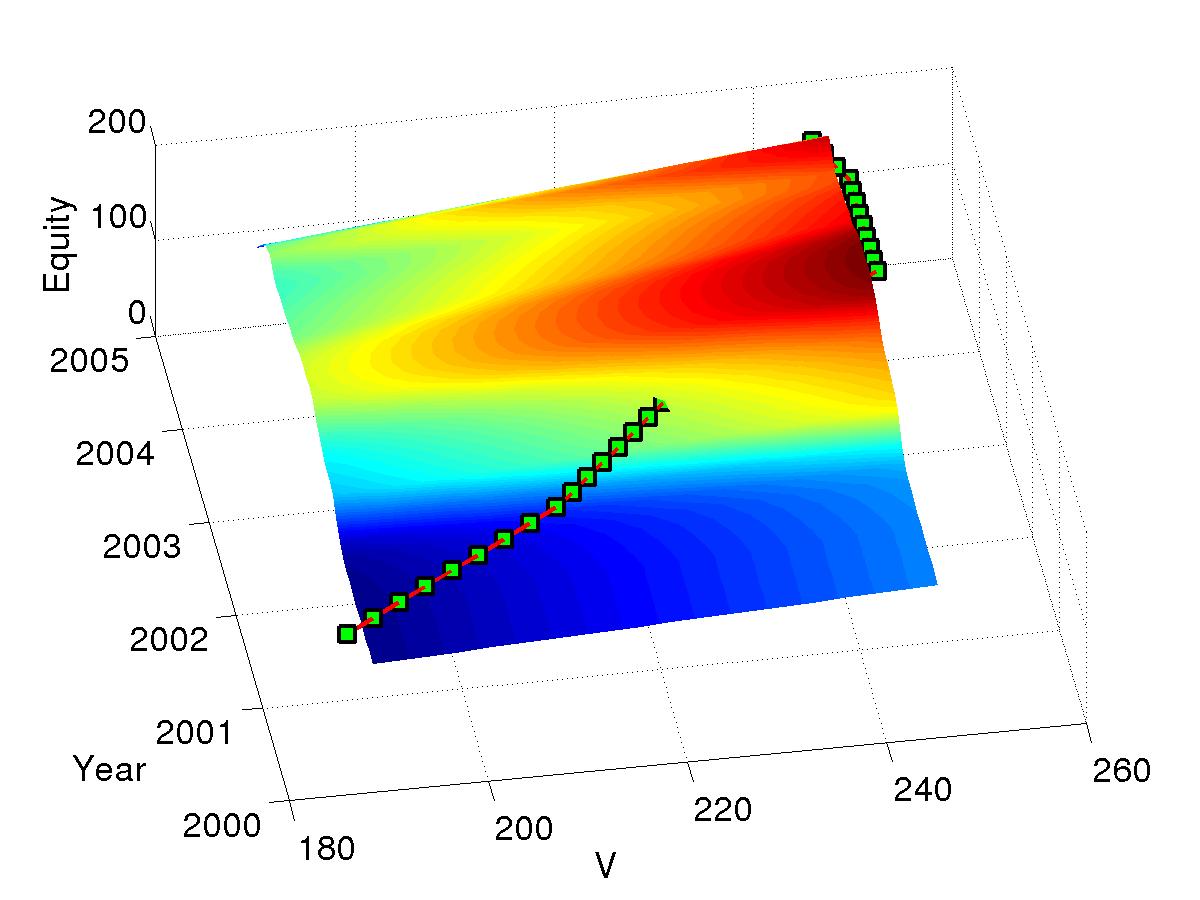}}
   \hskip 0.01\textwidth
   \subfigure[]{
   \label{FIG06d}
   \includegraphics[width=0.5\textwidth]{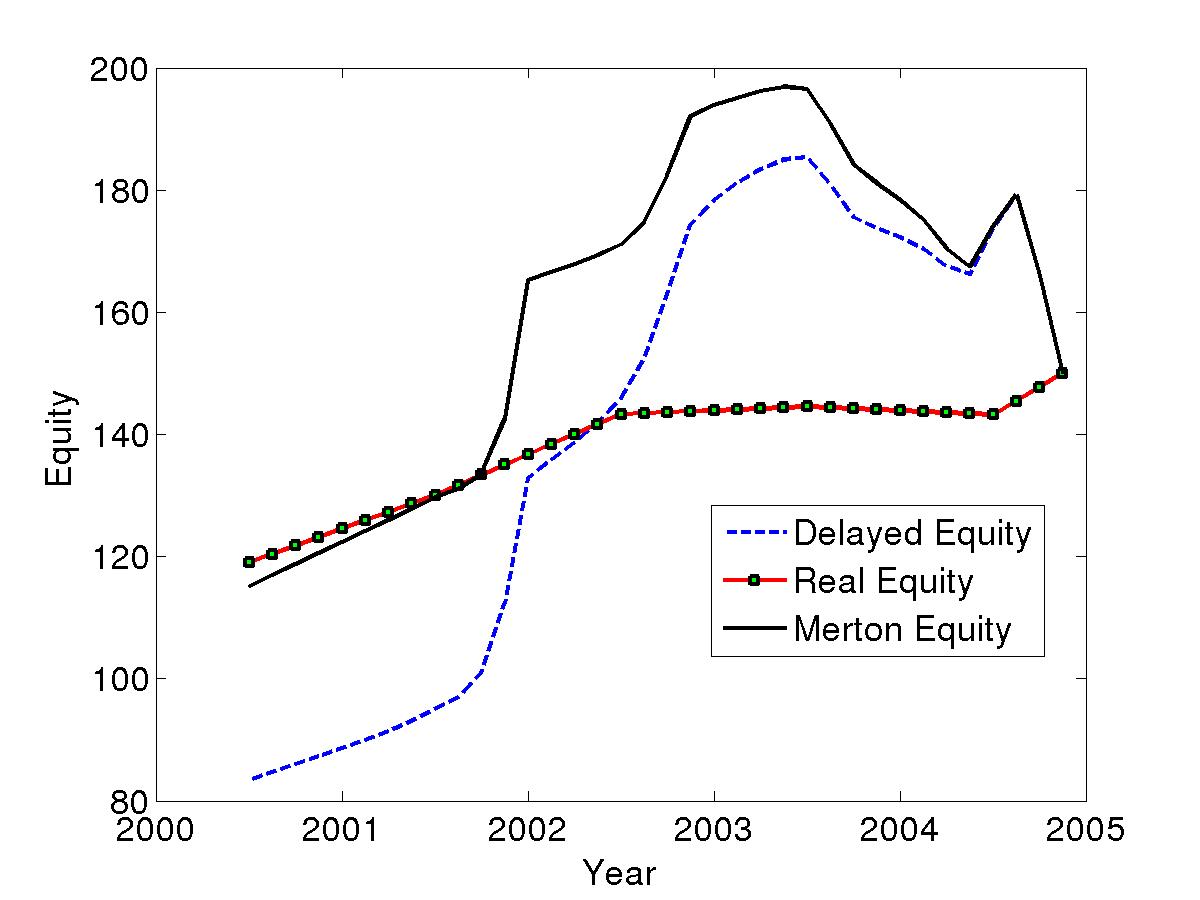}}
 \caption{ The graphs for  Firm $C_8$.  We plot at the left ((a) and (c) respectively)  the surface graphs of the numerical equity from our delayed model at $T=9.5$ and $T=5$.
 The corresponding  3 D graphs (green curves) of the real data of the firm equity value as a function of the time  and $V$ are also plotted in (a) and (c). 
 At the right ((b) and (d)), 2 D  graphs of the firm  equity  value as a function of time, corresponding to the surface graphs at  the left ((a) and (c)) respectively) are presented. Those 2 D graphs
 contain  the numerical equity from our delayed model, the numerical equity of the Merton model and the real data equity of the firm.}
 \label{FIG06}
\end{figure}

\begin{figure}[H]
 \subfigure[]{
   \label{FIG07a}
   \includegraphics[width=0.5\textwidth]{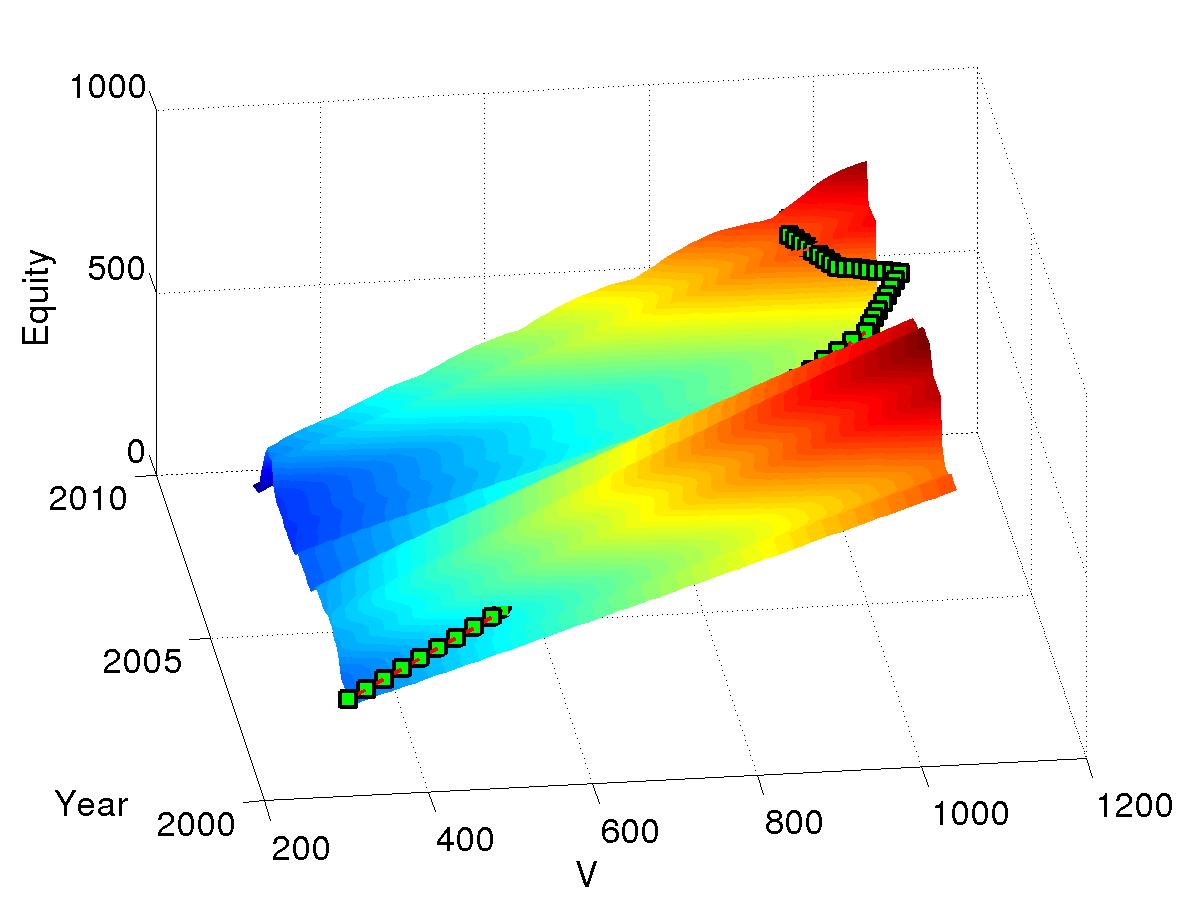}}
   \hskip 0.01\textwidth
   \subfigure[]{
   \label{FIG07b}
   \includegraphics[width=0.5\textwidth]{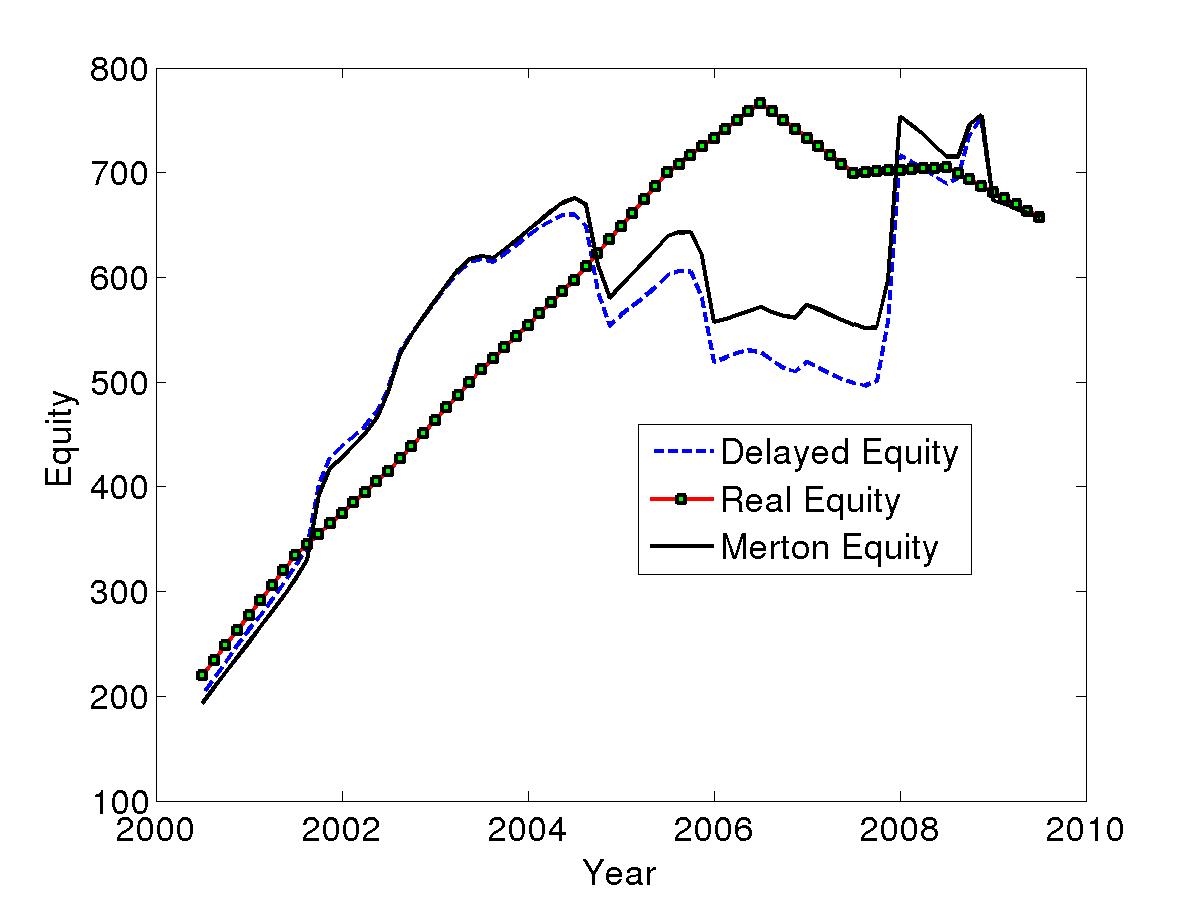}}
   \hskip 0.01\textwidth
   \subfigure[]{
   \label{FIG07c}
   \includegraphics[width=0.5\textwidth]{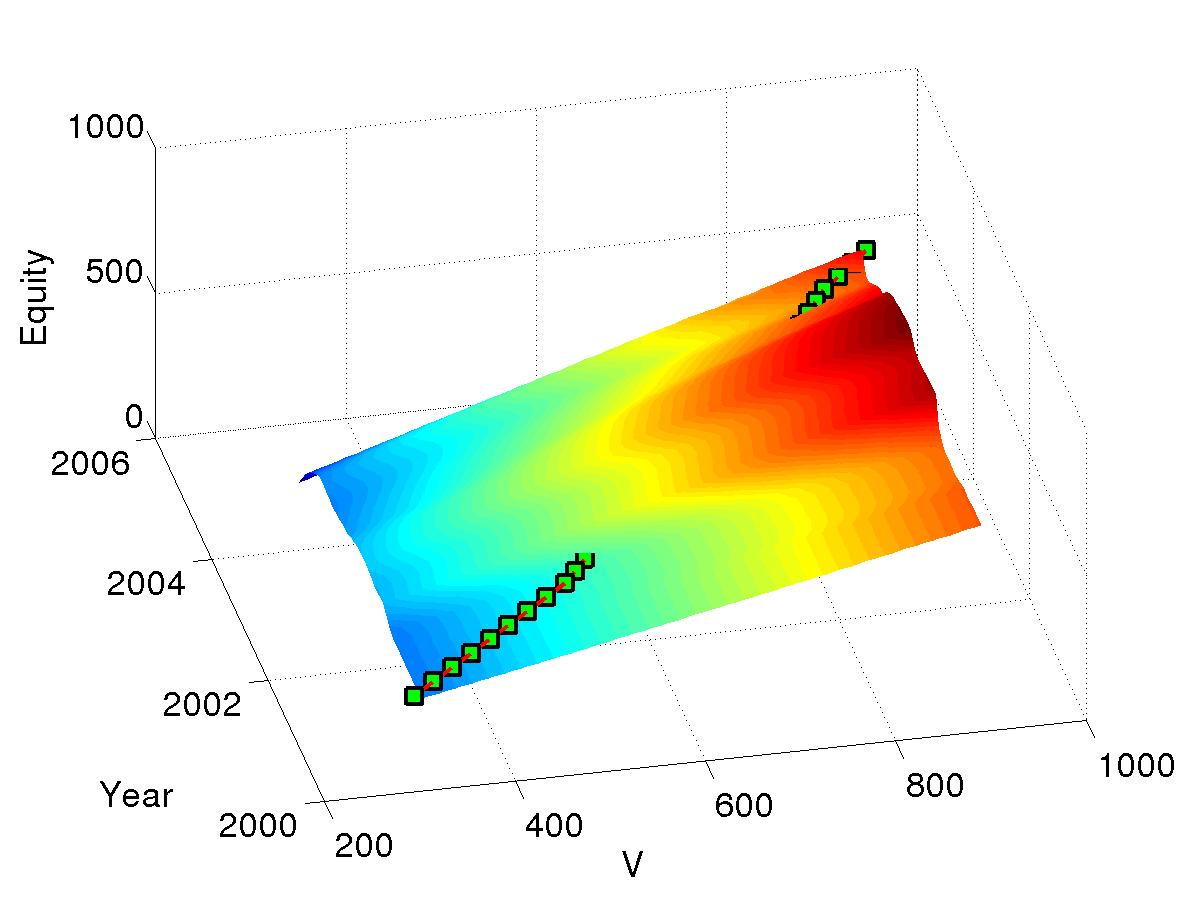}}
   \hskip 0.01\textwidth
   \subfigure[]{
   \label{FIG07d}
   \includegraphics[width=0.5\textwidth]{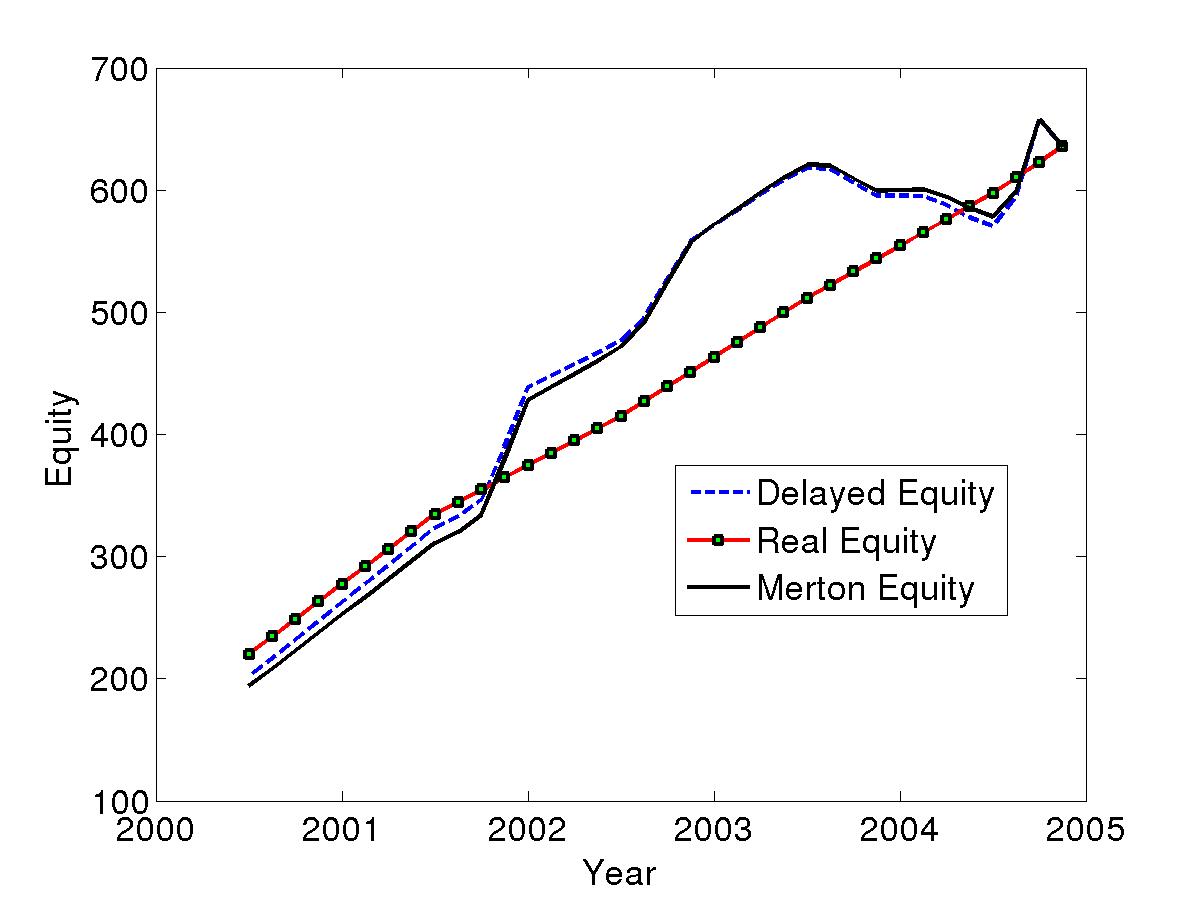}}
 \caption{ The graphs for  Firm $C_9$.  We plot at the left ((a) and (c) respectively)  the surface graphs of the numerical equity from our delayed model at $T=9.5$ and $T=5$.
 The corresponding  3 D graphs (green curves) of the real data of the firm equity value as a function of the time  and $V$ are also plotted in (a) and (c). 
 At the right ((b) and (d)), 2 D  graphs of the firm  equity  value as a function of time, corresponding to the surface graphs at  the left ((a) and (c)) respectively) are presented. Those 2 D graphs
 contain  the numerical equity from our delayed model, the numerical equity of the Merton model and the real data equity of the firm.}
 \label{FIG07}
\end{figure}

\begin{figure}[H]
 \subfigure[]{
   \label{FIG08a}
   \includegraphics[width=0.5\textwidth]{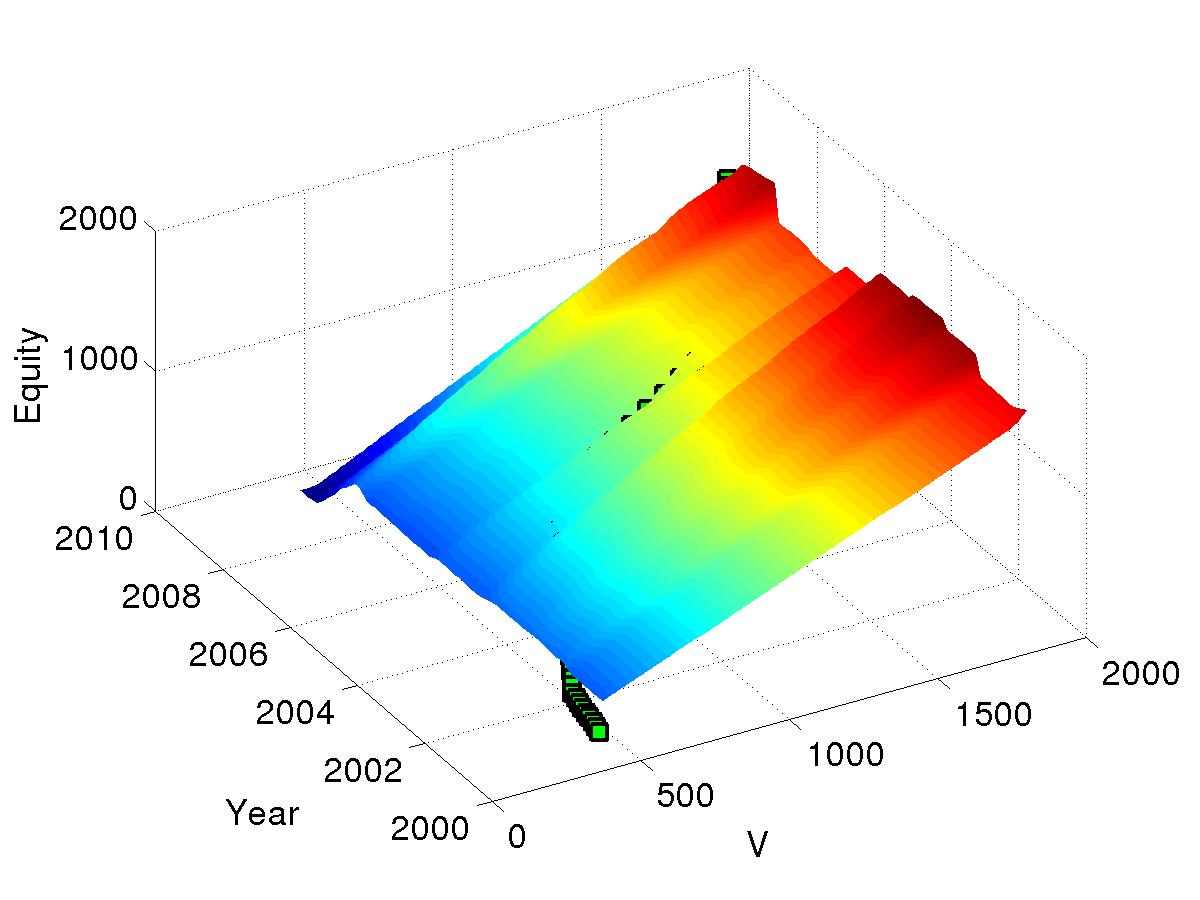}}
   \hskip 0.01\textwidth
   \subfigure[]{
   \label{FIG08b}
   \includegraphics[width=0.5\textwidth]{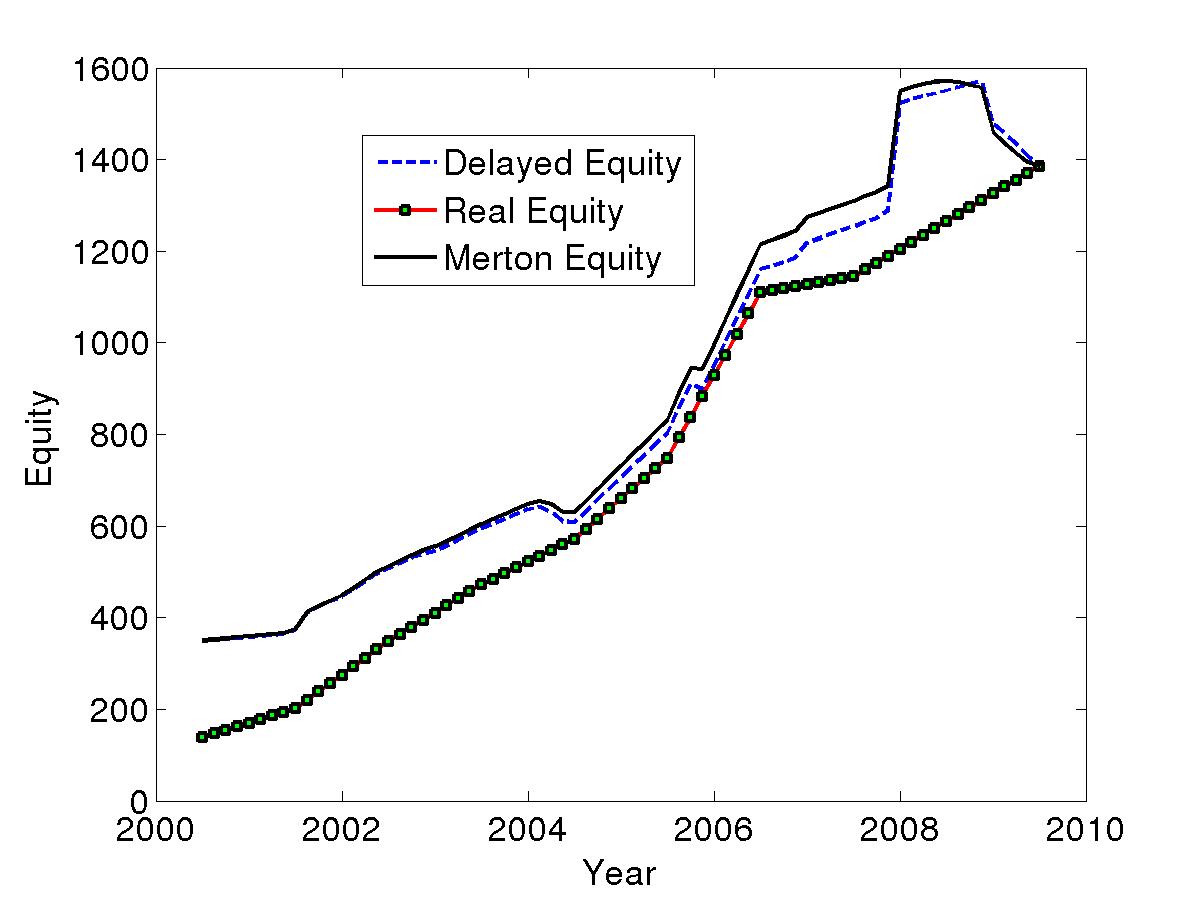}}
   \hskip 0.01\textwidth
   \subfigure[]{
   \label{FIG08c}
   \includegraphics[width=0.5\textwidth]{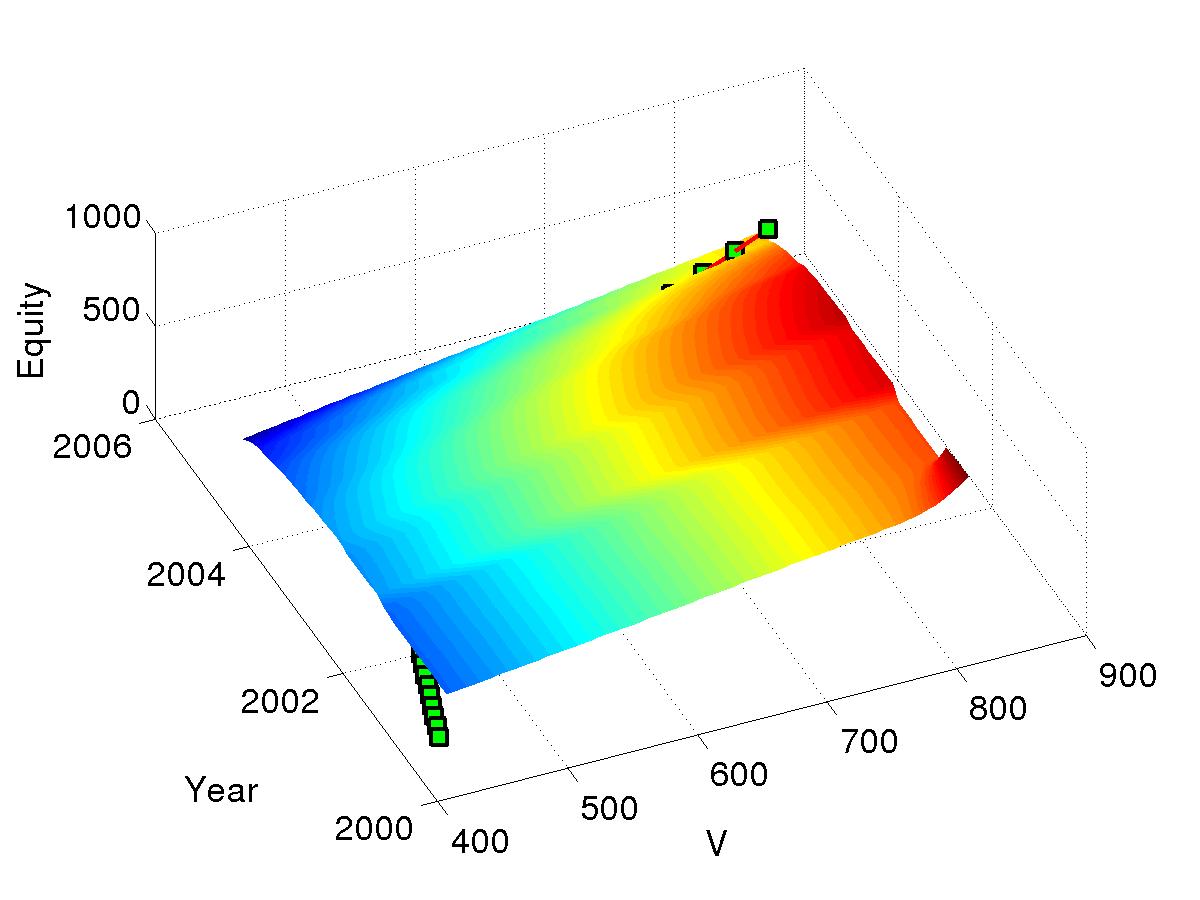}}
   \hskip 0.01\textwidth
   \subfigure[]{
   \label{FIG08d}
   \includegraphics[width=0.5\textwidth]{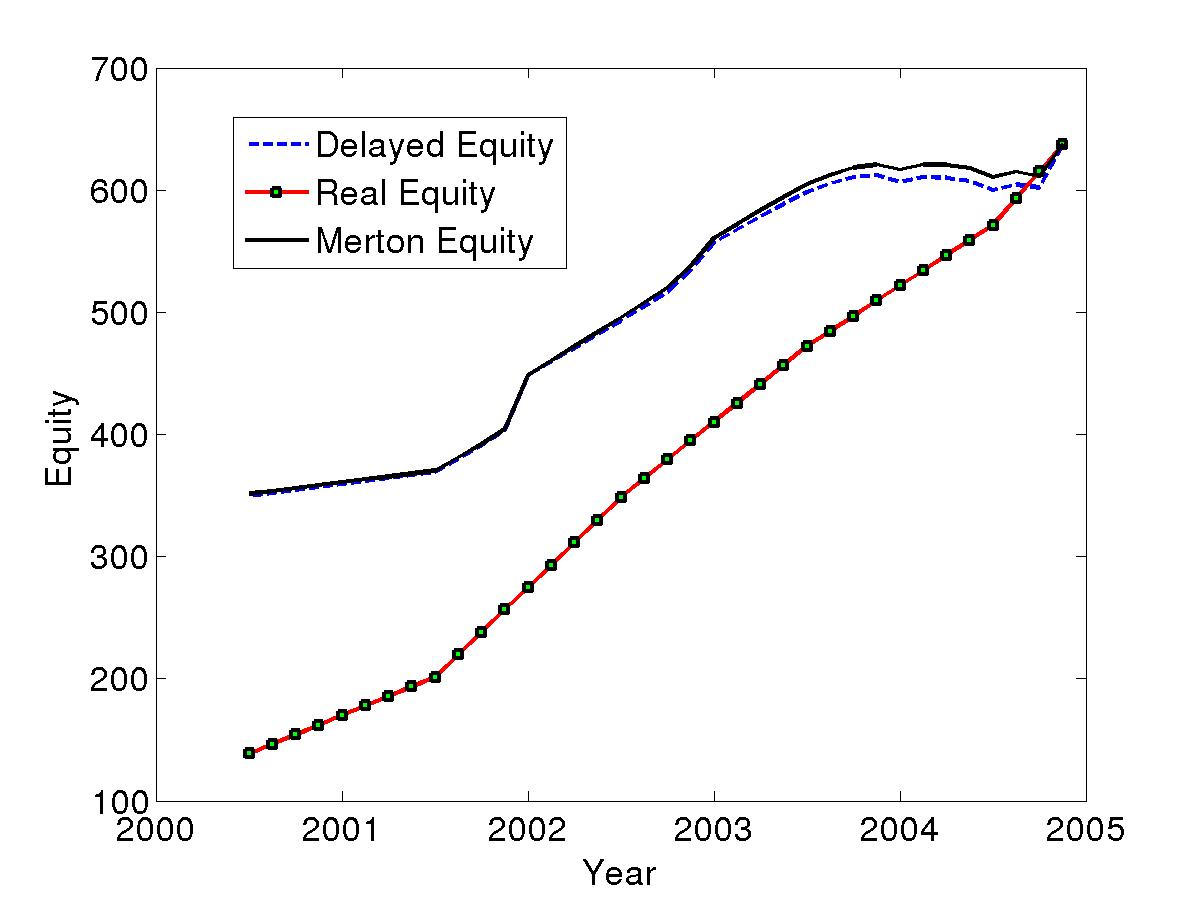}}
 \caption{ The graphs for  Firm $C_{10}$.  We plot at the left ((a) and (c) respectively)  the surface graphs of the numerical equity from our delayed model at $T=9.5$ and $T=5$.
 The corresponding  3 D graphs (green curves) of the real data of the firm equity value as a function of the time  and $V$ are also plotted in (a) and (c). 
 At the right ((b) and (d)), 2 D  graphs of the firm  equity  value as a function of time, corresponding to the surface graphs at  the left ((a) and (c)) respectively) are presented. Those 2 D graphs
 contain  the numerical equity from our delayed model, the numerical equity of the Merton model and the real data equity of the firm.}
 \label{FIG08}
\end{figure}

\begin{figure}[H]
 \subfigure[]{
   \label{FIG09a}
   \includegraphics[width=0.5\textwidth]{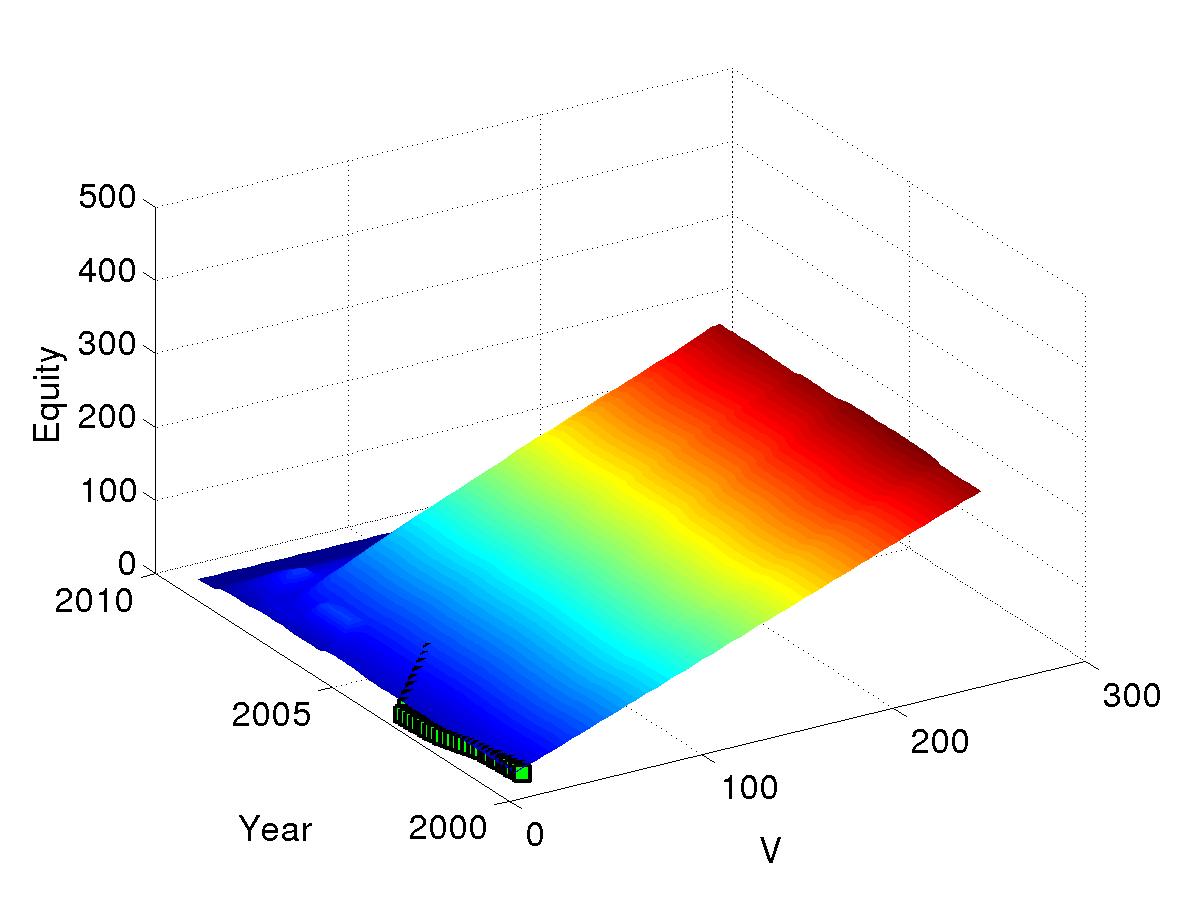}}
   \hskip 0.01\textwidth
   \subfigure[]{
   \label{FIG09b}
   \includegraphics[width=0.5\textwidth]{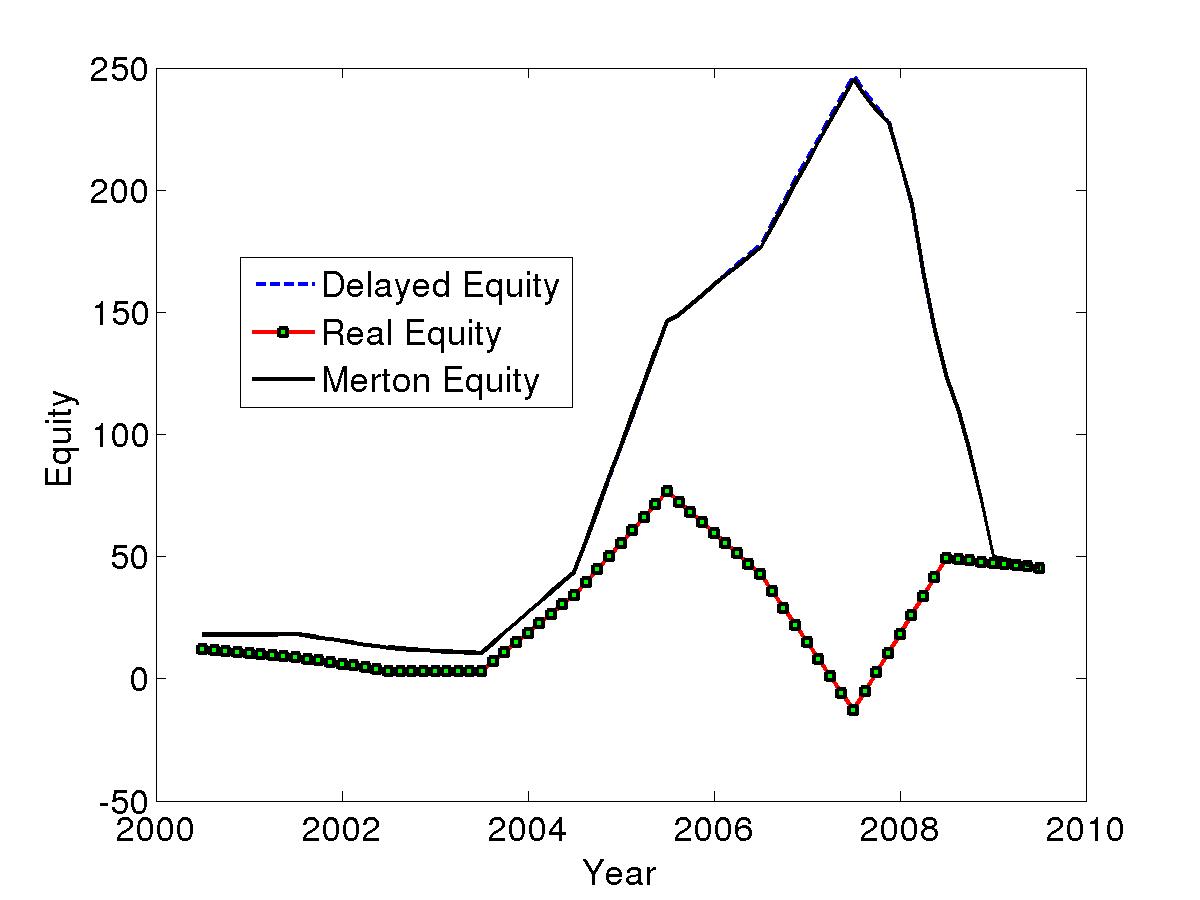}}
   \hskip 0.01\textwidth
   \subfigure[]{
   \label{FIG09c}
   \includegraphics[width=0.5\textwidth]{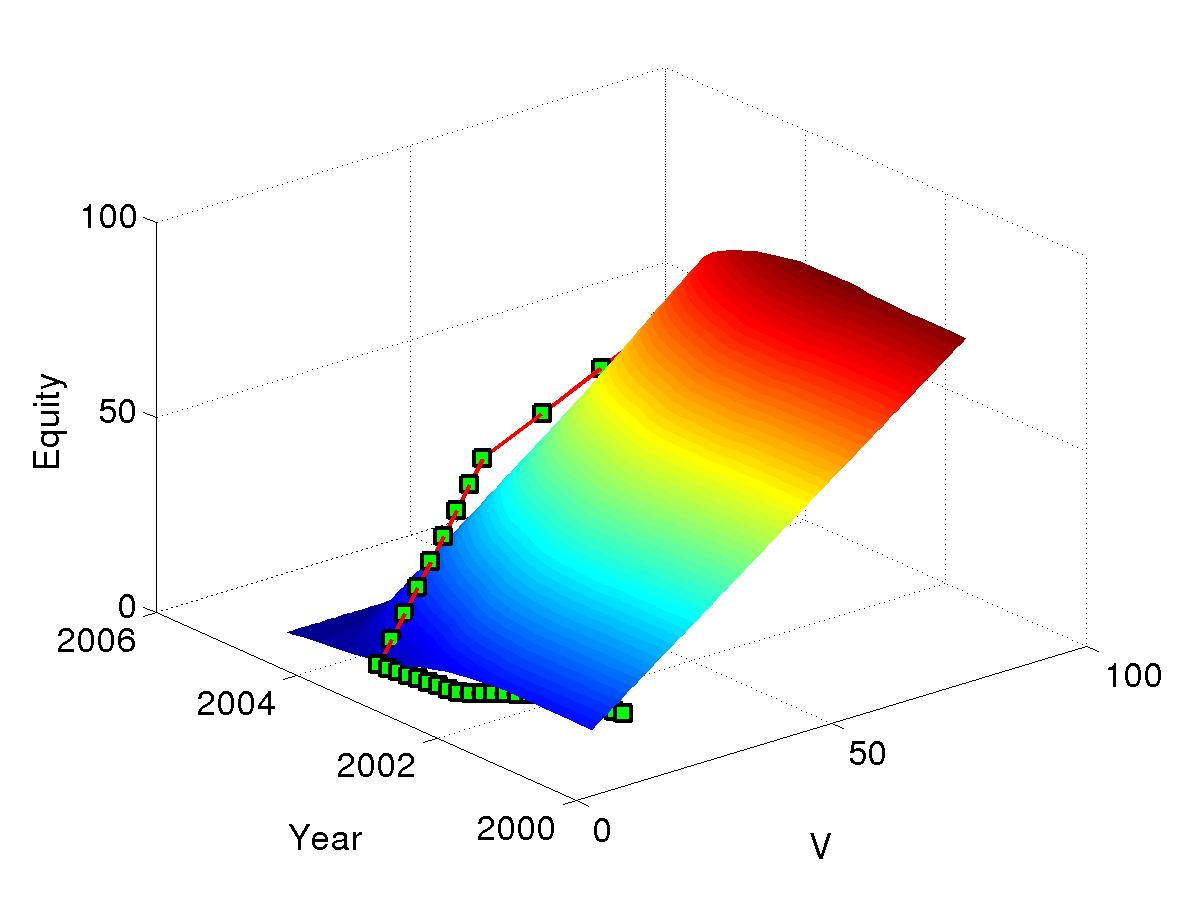}}
   \hskip 0.01\textwidth
   \subfigure[]{
   \label{FIG09d}
   \includegraphics[width=0.5\textwidth]{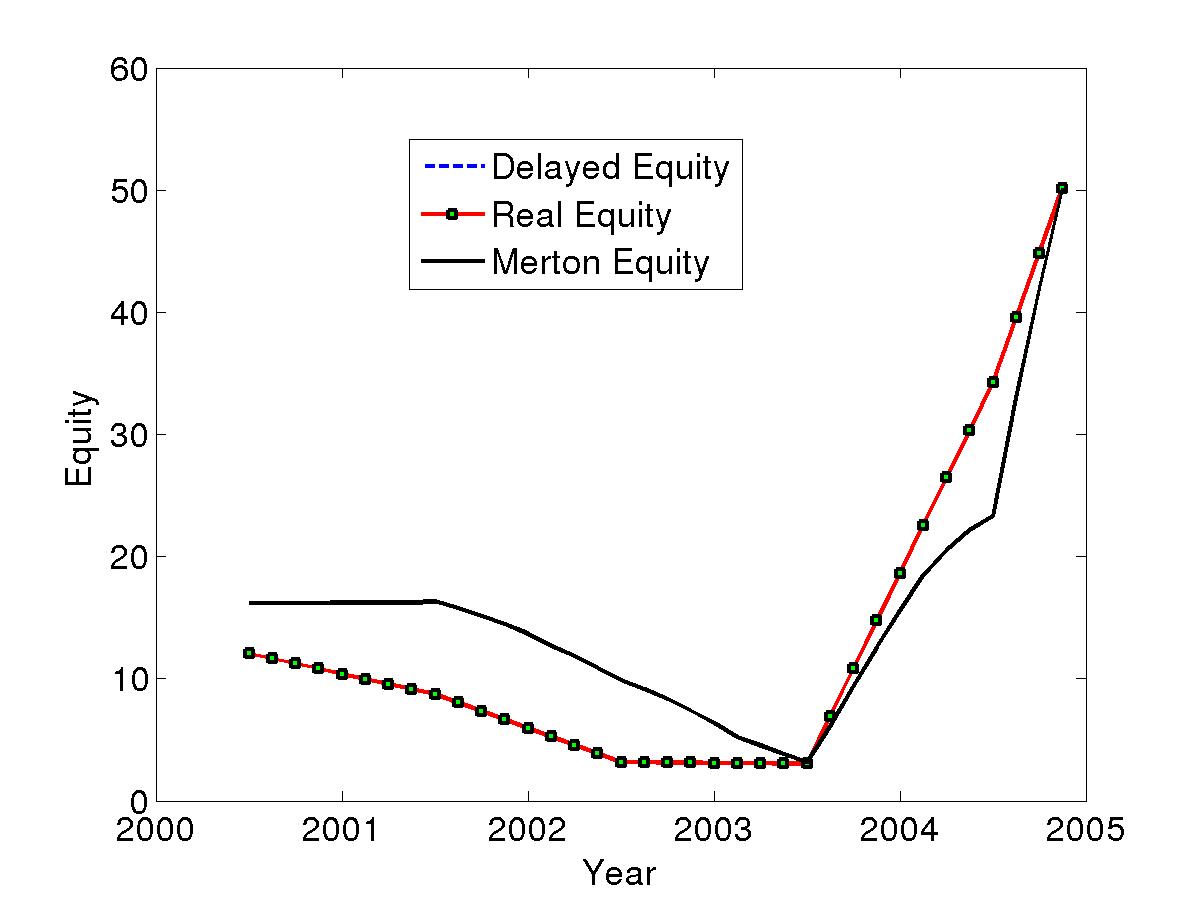}}
 \caption{ The graphs for  Firm $C_{11}$.  We plot at the left ((a) and (c) respectively)  the surface graphs of the numerical equity from our delayed model at $T=9.5$ and $T=5$.
 The corresponding  3 D graphs (green curves) of the real data of the firm equity value as a function of the time  and $V$ are also plotted in (a) and (c). 
 At the right ((b) and (d)), 2 D  graphs of the firm  equity  value as a function of time, corresponding to the surface graphs at  the left ((a) and (c)) respectively) are presented. Those 2 D graphs
 contain  the numerical equity from our delayed model, the numerical equity of the Merton model and the real data equity of the firm.}
 \label{FIG09}
\end{figure}

We used  the following frims:
\begin{itemize}
 \item [$C_5$] {\tt Magna International Inc} (\figref{FIG03})
 \item [$C_6$]  {\tt First Citizens Bancshares Inc NC} (\figref{FIG04})
\item [$C_7$]  {\tt Coca-Cola CO} (\figref{FIG05})
\item [$C_8$] {\tt One Liberty Properties INC } (\figref{FIG06})
\item [$C_9$] {\tt Cisco Systems INC} (\figref{FIG07})
\item [$C_{10}$] {\tt C B S Corp NEW } (\figref{FIG08})
\item [$C_{11}$] {\tt  Nam Tai Electronics INC} (\figref{FIG09})
\end{itemize}

  Here again,  the time origin corresponds to the year (2000+1/2), the data before
are memory data and we want to predict the data after (2000+1/2).
In the legends of all of our graphs we use the following notation
\begin{itemize}
  \item ``Delayed Equity'' is for the numerical equity value from our nonlinear  delayed model.
  \item  ``Real Equity'' is for the real equity value of the corporate.
  \item ``Merton Equity'' is for the numerical equity value from Merton model.
\end{itemize}
In our  surface graphs of  the  numerical equity value,  we plot only the part  where the variable  $V$ is between the minimun and the maximum values of  our real market value $V$.
In all simulations with our delayed model,  we take $L=9.5$.
In all  graphs, the function $g$ (volatility in delayed model) is the quadratic interpolation 
of the  standard deviation of daily returns $\sigma$ in the memory part while the volatility in the Merton model is just the mean of the memory part.

For each firm, we plot at the left hand size  both the surface graphs of the numerical equity value from our delayed model at $T=9.5$ and $T=5$. In those 3D surface graphs,  
we also plot the corresponding  3 D graphs (green curves) of the real data of the firm equity value as a function of the time (year) and $V$. 
At the right hand size, we plot in 2 D  the firm  equity  value as a function of time (year), corresponding to the surface  graphs at the left hand size.
Those 2 D equity graphs contain the numerical equity value from our delayed model, the numerical equity value of the Merton model and the real data equity value of the firm.

In our simulations, for a given $T$,  the promised debt $B$ is just the real debt value of the firm at time $T$. 

For firm $C_5$ in \figref{FIG03}, we can observe that both the delayed model and Merton model fit well the real market equity value of the firm. The accuracy of the two methods varies within some time interval
as we can observe in \figref{FIG03b} and \figref{FIG03d}. 

For firm $C_6$ in \figref{FIG04}, comparing to firm $C_5$ the two models fit less.  In a wide time interval  in  \figref{FIG04b} and \figref{FIG04d}, the delayed  model is more close to the real market equity of the firm.
We can also observe a good early fit in the Merton model.

For firm $C_7$ in \figref{FIG05}, comparing to firm $C_5$ the two models fit less. But for the maturity date $T=9.5$ in   \figref{FIG05} the fitting is relatively good for the two  models.
The accuracy of the two methods varies within some time interval as we can observe in \figref{FIG05b} and \figref{FIG05d}. 

For firm $C_8$ in \figref{FIG06}, the fitting is relatively bad for the two models. However we can observe in \figref{FIG06b} and \figref{FIG06d} the good early fit in the Merton model, and that 
in the wide time interval the delayed model is more close to the real market equity of the firm than the Merton model.

For firm $C_9$ in \figref{FIG07}, the fitting is relatively good for the two models in the early time interval and become relatively bad just after.

For firm $C_{10}$ in \figref{FIG08}, the fitting is relatively good for the two models for the maturity date $T=9.5$ in \figref{FIG08b} at the middle time interval  and  bad for the maturity date $T=5$ in \figref{FIG08d}.

For firm $C_{11}$ in \figref{FIG09}, the fitting is relatively good for the two models in the early time interval but become bad just after. The two models are confused.


\section{Conclusion}
\label{conclusion}
 In this paper,  numerical techniques to solve delayed nonlinear model  for  pricing  corporate liabilities  are provided.
The numerical technique to solve  the  RPDEs modeling debt and equity value combines the finite difference--finite volume methods 
(discretization  respect to the firm value $V$) and  an exponential integrator (discretization respect to the time $t$).
 The matrix exponential functions are computed efficiently using Krylov subspace technique.

Using  financial data from some firms, we compare  numerical solutions from both our nonlinear model and classical  Merton  to the real firm's data.
This comparaison  shows that our nonlinear model behaves  very well.  We  conclude that in corporate  finance  the past dependence of
the firm value process may be an important feature and therefore  should not be ignored.

\section*{ACKNOWLEDGEMENTS}
We thank  Dr. David Rakwoski from   College of Business, Southern Illinois University
 for  finding data for the simulations. 
 Antoine  Tambue  was funded by the Research Council of
Norway (grant number 190761/S60).

\newpage

\newpage

\newpage

\newpage
\newpage
\newpage
\newpage
\newpage
\newpage 



\end{document}